# ON A REFINED STARK CONJECTURE
# FOR FUNCTION FIELDS

Cristian D. Popescu

ABSTRACT. We prove that a refinement of Stark's Conjecture formulated by Rubin in [14] is true up to primes dividing the order of the Galois group, for finite, abelian extensions of function fields over finite fields. We also show that in the case of constant field extensions a statement stronger than Rubin's holds true.

## 0. Introduction

Let $L_{K/k}(s,\chi)$ be the Artin $L$–function associated to a finite, Galois extensions of global fields $K/k$, and a character $\chi \in \widehat{G(K/k)}$. In the 1970's and early 1980's, Stark [15] developed a conjecture concerning the leading coefficient of the Taylor expansion of $L_{K/k}(s,\chi)$, at $s = 0$. The classical formulae for the case of a Dedekind Zeta–function provide good hints: if the order of vanishing is $r$, then the coefficient in question should be a rational number multiplied by a regulator, obtained as a determinant of an $r \times r$ matrix involving character values and logarithms of absolute values of units in $K$.

The examples however suggest that at least in the case of *abelian L*-functions, the denominator of the above mentioned rational number could be specified. In the last paper of [15], Stark proposed a refined conjecture, which makes this denominator specific, in the case of *abelian L*-functions of order of vanishing $r = 1$ at $s = 0$. This is what Tate (see [16]) calls a conjecture "over **Z**" and he asks for such a refined statement for abelian $L$–functions of any order of vanishing $r$. In [14] Rubin formulates a refined version of Stark's conjecture, for abelian $L$–functions of number fields, of any order of vanishing $r$ at $s = 0$. Rubin's statement can be easily extended to the more general case of abelian extensions of global fields of any characteristic.

In this paper we study Rubin's conjecture in the case of finite, abelian extensions of global fields of characteristic $p > 0$ (i.e. function fields). The case $r = 1$ in this situation was independently solved by Deligne (see [17, Chpt.V]) and Hayes (see [6]), with methods relying on $\ell$–adic homology of 1–motives and rank 1 sign–normalized Drinfeld modules respectively. We adopt methods similar to Deligne's in order to treat the arbitrary order of vanishing case.

In §1 we set the notation, define the objects involved and some of their properties, and state the conjecture. §2 is concerned with the "dynamics" of Rubin's conjecture. In particular we show that, in certain cases, one can derive the conjecture at any level $r$ from a stronger form of the conjecture at level $r = 0$ (see Corollary

Research at MSRI is supported in part by NSF grant DMS–9022140

Typeset by $\mathcal{AMS}$-TEX

2.2). This fact turns out to be of crucial importance in §§3 and 4. In §3 we prove Rubin's conjecture, up to primes dividing the order of the Galois group $G(K/k)$, for any abelian extension $K/k$ of function fields (see Theorems 3.1.1 and 3.2.1). In §4 we prove that a stronger statement than Rubin's holds true for constant field extensions of function fields (see Theorems 4.2.9 and 4.3.1).

We use these results in [13] in order to formulate and prove Gras–type Conjectures for global fields of characteristic $p > 0$.

The main technique employed in this is the $\ell$–adic homological (or, equivalently, the $\ell$-adic étale cohomological) ($\ell \neq p$) interpretation of $L$–functions, due to Grothendieck. However, in order to deal with the $p$–part of the conjecture, we had to rely on the $p$–adic homology groups, which are known not to give a good theory of $L$–functions. Therefore, we needed a "bridge" between the $\ell$–adic étale and the $p$–adic étale cohomology theories, and that is provided by the crystalline cohomology theory. The fact that the Frobenius actions on the $\ell$-adic étale and the crystalline cohomology groups have the same characteristic polynomials has been known for quite some time (see [8]). However, we needed the slightly stronger statement that the characteristic polynomials of the Frobenius actions on the $\chi$–components of the $\ell$–adic étale and crystalline cohomology groups are the same, for every character $\chi$ of $G(K/k)$. Due to lack of a good reference, we provide a proof of this fact (see Proposition A.1) along with a proof of Theorem 1.7.4.1, in the Appendix.

**Acknowledgements.** It is a pleasure to thank Karl Rubin for his invaluable help and support. Our approach to this problem has been greatly influenced by his published and unpublished thoughts and by John Tate's book [17].

## 1. The objects involved and their properties

1.1. GENERAL NOTATION

Let us fix a function field $k$ of characteristic $p > 0$ and of transcendental degree 1 over a finite field. Let $K$ be a finite, abelian extension of $k$, of Galois group $G = G(K/k)$, and let $g = |G|$. We denote by $\mathbf{F}_q$ and $\mathbf{F}_{q^\nu}$ the exact fields of constants of $k$ and $K$ respectively, where $q$ is a power of $p$, and $\nu$ is a positive integer. For primes $v$ in $k$ and $w$ in $K$, such that $w|v$, we denote by $\mathbf{F}_q(v)$ and by $\mathbf{F}_{q^\nu}(w)$ their corresponding residue fields, and by $d_v$ and $d_w$ their degrees over $\mathbf{F}_q$ and $\mathbf{F}_{q^\nu}$ respectively ($d_v = [\mathbf{F}_q(v) : \mathbf{F}_q]$ and $d_w = [\mathbf{F}_{q^\nu}(w) : \mathbf{F}_{q^\nu}]$). Let $\mathrm{N}v = |\mathbf{F}_q(v)|$, $\mathrm{N}w = |\mathbf{F}_{q^\nu}(w)|$, and let $G_v$ be the decomposition group associated to $v$ in $K/k$. The normalized absolute value $|\cdot|_w : K^\times \longmapsto q^{\mathbf{Z}}$ associated to $w$ is defined by: $|\alpha|_w = (\mathrm{N}w)^{-ord_w(\alpha)}$, for every $\alpha \in K$.

Let $S$ and $T$ be two finite, nonempty and disjoint sets of primes in $k$, $S$ containing all primes which ramify in $K/k$. Let us define:
• $S_K = \{$primes in $K$ lying above primes in $S\}$
• $T_K = \{$primes in $K$ lying above primes in $T\}$

- $O_S = \{\alpha \in K : |\alpha|_w \leq 1, \text{ for all } w \notin S_K\}$
- $U_S = O_S^\times$ (the group of $S$-units in $K$).
- $U_{S,T} = \{\alpha \in U_S : \alpha \equiv 1 \mod w, \text{ for all } w \in S_K\}$ (This is a free $\mathbf{Z}$–module.)
- $A_S = $ the ideal class group of $O_S$
- $A_{S,T} = \dfrac{\{\text{Fractional ideals of } O_S, \text{ prime to } T_K\}}{\{fO_S : f \equiv 1 \mod w, \forall w \in T_K\}}$. $A_{S,T}$ is a finite ideal class group, isomorphic via the Artin reciprocity map with the Galois group of the maximal abelian extension of $K$, of conductor dividing $\prod_{w \in T_K} w$, in which all primes in $S_K$ split completely.
- $X_S = \{\sum_{w \in S_K} a_w \cdot w : a_w \in \mathbf{Z}, \sum_{w \in S_K} a_w = 0\}$
- $\lambda_{S,T} : U_{S,T} \longrightarrow \mathbf{R} \otimes X_S$, is the $G$–morphism defined by

$$\lambda_{S,T}(\alpha) = \sum_{w \in S_K} -\log(|\alpha|_w) \cdot w,$$

for any $\alpha \in U_{S,T}$. This induces an $\mathbf{R}[G]$–isomorphism $\mathbf{R} \otimes U_{S,T} \xrightarrow{\sim} \mathbf{R} \otimes X_S$.
- $R_{S,T}$ is the absolute value of the determinant of $\lambda_{S,T}$, with respect to $\mathbf{Z}$–bases of $U_{S,T}$ and $X_S$.

There is an exact sequence of $\mathbf{Z}[G]$–modules:

$$0 \longrightarrow U_{S,T} \longrightarrow U_S \longrightarrow \bigoplus_{w \in T_K} \mathbf{F}_{q^\nu}(w)^\times \longrightarrow A_{S,T} \longrightarrow A_S \longrightarrow 0, \qquad (1)$$

which, together with the usual class–number formula, implies that, if

$$\zeta_{S,T}(s) = \prod_{w \notin S_K} (1 - \mathrm{N}w^{-s})^{-1} \cdot \prod_{w \in T_K} (1 - \mathrm{N}w^{1-s})$$

is the $(S,T)$–zeta function associated to $K$, then:

$$\mathrm{ord}_{s=0}\zeta_{S,T}(s) = |S_K| - 1 \quad , \quad \lim_{s \to 0} s^{1-|S_K|}\zeta_{S,T}(s) = -|A_{S,T}| \cdot R_{S,T}. \qquad (2)$$

(see [4] for (1) and (2) above)

1.2. THE $L$–FUNCTIONS

Let $\widehat{G}$ be the group of irreducible, complex valued characters of $G$. For any $\chi \in \widehat{G}$, $L(s,\chi)$ will denote the Artin $L$–function associated to $\chi$, and $e_\chi = 1/g \sum_{\sigma \in G} \chi(\sigma) \cdot \sigma^{-1} \in \mathbf{C}[G]$. The Stickelberger function is defined by:

$$\Theta(s) = \sum_{\chi \in \widehat{G}} L(s,\chi) \cdot e_{\chi^{-1}},$$

and it can be thought of as a complex meromorphic function, with values in $\mathbf{C}[G]$.

Rather than working with the $L$ and Stickelberger functions as defined above, Rubin [14] works with the $(S,T)$–modified objects, defined by:

$$L_{S,T}(s,\chi) = \prod_{v\in T}(1-\chi(\sigma_v)\cdot Nv^{1-s}) \cdot \prod_{v\notin S}(1-\chi(\sigma_v)\cdot Nv^{-s})^{-1}$$

$$\Theta_{S,T}(s) = \sum_{\chi\in\widehat{G}} L_{S,T}(s,\chi)\cdot e_{\chi^{-1}},$$

where $\sigma_v$ is the Frobenius morphism associated to (the unramified) $v$ in $G$. Due to the conditions imposed on $S$ and $T$, $L_{S,T}(s,\chi)$ and $\Theta_{S,T}(s)$ are holomorphic on the whole complex plane. For an integer $r \geq 0$, such that $s^{-r}\Theta_{S,T}(s)$ is holomorphic at $s=0$, let

$$\Theta_{S,T}^{(r)}(0) = \lim_{s\to 0} s^{-r}\Theta_{S,T}(s) = \sum_{\chi\in\widehat{G}} \lim_{s\to 0} s^{-r} L_{S,T}(s,\chi)\cdot e_{\chi^{-1}} \in \mathbf{C}[G].$$

For $S$ and $T$ fixed, and $\chi \in \widehat{G}$, let $r_\chi = \text{ord}_{s=0} L_{S,T}(s,\chi)$. These numbers do not depend on $T$, but they depend on $S$ in the following way:

$$r_\chi = \begin{cases} \#\{v\in S : \chi|_{G_v} = 1_{G_v}\}, & \text{if } \chi \neq 1_G \\ \#S - 1, & \text{if } \chi = 1_G, \end{cases} \qquad (3)$$

or equivalently (and more elegantly)

$$r_\chi = \langle \chi, \chi_{X_S}\rangle, \text{ for all } \chi \in \widehat{G}, \qquad (4)$$

where $\chi_{X_S}$ is the character associated to the Galois representation $\rho_{X_S}$ given by the $\mathbf{C}[G]$–module $\mathbf{C}\otimes X_S$, and $\langle\,,\,\rangle$ is the usual inner product on the space of characters $\widehat{G}$. (see [17] for (3) and (4) above.)

1.3. GROUP RINGS AND $G$–MODULES

For a subring $R$ of $\mathbf{C}$, and a $\mathbf{Z}[G]$–module $M$, $R[G]$ denotes the group ring with coefficients in $R$, and $RM$ denotes the $R[G]$–module $R\otimes_{\mathbf{Z}} M$. For $\chi \in \widehat{G}$, $R[\chi]$ denotes the minimal ring extension of $R$, containing the values of $\chi$.

If $L$ is a field of characteristic $0$, then $\widehat{G}(L)$ will denote the set of characters associated to $L$–irreducible representations of $G$. (In particular, if $L = \mathbf{C}$, then $\widehat{G}(\mathbf{C}) = \widehat{G}$.) If $\bar{L}$ is an algebraic closure of $L$, then $G(\bar{L}/L)$ acts canonically on $\widehat{G}(\bar{L})$ and $\widehat{G}(L)$ can be viewed as the set of orbits with respect to this action. For $\psi \in \widehat{G}(L)$ and $\chi \in \widehat{G}(\bar{L})$, we write $\chi|\psi$ if $\chi$ belongs to the orbit represented by $\psi$. If $S$ and $T$ are fixed, (3) above implies that $r_\chi = r_{\chi^\tau}$, for any $\chi \in \widehat{G}(\bar{L})$, and any $\tau \in G(\bar{L}/L)$. We can therefore define $r_\psi$ as $r_\psi = r_\chi$, for any $\psi \in \widehat{G}(L)$ and any $\chi \in \widehat{G}(\bar{L})$, such that $\chi|\psi$.

If $r$ is a positive integer and $S$ is fixed, let $\widehat{G}(L, r) = \{\psi \in \widehat{G}(L) : r_\psi = r\}$. If $R$ is a ring whose field of fractions is $L$, and $M$ is an $R[G]$– module, let

$$M_{r,S} = \{m \in M : e_\psi \cdot m = 0, \text{ for all } \psi \in \widehat{G}(L) \setminus \widehat{G}(L, r)\},$$

where $e_\psi = 1/g \sum_{\sigma \in G} \psi(\sigma) \cdot \sigma^{-1}$. This is obviously an $R[G]$–submodule of $M$.

**Remark 1.** Let $R$ be a Dedekind domain containing $\mathbf{Z}[1/g]$, and let $L$ be its field of fractions. Then one has a decomposition

$$R[G] = \bigoplus_{\psi \in \widehat{G}(L)} D_\psi,$$

where $D_\psi = R[G] \cdot e_\psi$. The $D_\psi$'s are finite extensions of $R$ (and therefore Dedekind domains themselves), and $D_\psi \xrightarrow{\sim} R[\chi]$ via the map $R[G] \xrightarrow{\chi} \mathbf{R}[\chi]$, defined as $\chi(\sum_{\sigma \in G} a_\sigma \cdot \sigma) = \sum_{\sigma \in G} a_\sigma \chi(\sigma)$, for any $\chi \in \widehat{G}(\bar{L})$, $\chi|\psi$. This implies that, if $M$ is an $R[G]$–module, then one can decompose it into its $\psi$–components

$$M = \bigoplus_{\psi \in \widehat{G}(L)} M^\psi,$$

where $M^\psi = M \otimes_R D_\psi = \{m \in M : e_\chi \cdot m = 0 \text{ in } \bar{L}M, \text{ if } \chi \nmid \psi\}$. For a finitely generated $R[G]$–module $M$, the following are therefore equivalent:

$(i)$    $M$ is a projective $R[G]$–module.
$(ii)$   $M^\psi$ is a projective $D_\psi$–module, for any $\psi$.
$(iii)$   $M^\psi$ has no $D_\psi$–torsion, for any $\psi$.
$(iv)$   $M$ has no $R$–torsion.

1.4. FITTING IDEALS

If $A$ is a commutative, Noetherian ring, and $M$ is a finitely generated $A$–module, then $\text{Fitt}_A(M)$ denotes its Fitting ideal. If

$$A^n \xrightarrow{\phi} A^m \longrightarrow M \longrightarrow 0$$

is a finite presentation of $M$, one can consider the composition of $A$–morphisms

$$\bigwedge^m A^n \xrightarrow{\bigwedge^m \phi} \bigwedge^m A^m \xrightarrow{\det} A.$$

By definition, we have $\operatorname{Fitt}_A(M) = \operatorname{Im}(\det \circ \overset{m}{\wedge} \phi)$.

We will use the following properties of Fitting ideals:

(a) If $M$ is a cyclic $A$–module then $\operatorname{Fitt}_A(M) = \operatorname{Ann}_A(M)$.

(b) If $A \overset{f}{\to} B$ is a morphism of Noetherian rings, and $M$ is a finitely generated $A$–module, then $\operatorname{Fitt}_B(M \underset{A}{\bigotimes} B) = f(\operatorname{Fitt}_A(M))B$.

(c) If $R$ is a Dedekind domain, $M$ is a finitely generated $R[G]$–module, and $a \in R[G]$, then
$$a \in \operatorname{Fitt}_{R[G]}(M) \iff a \in \operatorname{Fitt}_{R_\lambda[G]}(M \underset{R}{\bigotimes} R_\lambda), \quad \forall \lambda,$$
where $\lambda$ runs through the set of prime ideals of $R$.

(d) If $A$ is a finite, direct sum of Noetherian rings $A = \bigoplus_i A_i$, $N = \bigoplus_i N_i$ is an $A$–module, with $N_i \overset{\sim}{\to} A_i^{n_i}$, for some positive integers $n_i$, and
$$N = \bigoplus_i N_i \xrightarrow{f = \bigoplus_i f_i} N = \bigoplus_i N_i \longrightarrow M \longrightarrow 0$$
is an exact sequence of $A$–modules, then $\det(f) \overset{\text{def}}{=} \sum_i \det_{A_i}(f_i) \in \operatorname{Fitt}_A(M)$. (Here $\det_{A_i}(f_i)$ is the determinant of $f_i$ with respect to an $A_i$ basis of $N_i$.)

Let $0 \longrightarrow M' \longrightarrow M \longrightarrow M'' \longrightarrow 0$ be an exact sequence of finitely generated $A$–modules. Then:

(e) $\operatorname{Fitt}_A(M') \cdot \operatorname{Fitt}_A(M'') \subseteq \operatorname{Fitt}_A(M) \subseteq \operatorname{Fitt}_A(M'')$.

(f) If $A$ is a Dedekind domain and $M$, $M'$ and $M''$ are finite, then
$$\operatorname{Fitt}_A(M) = \operatorname{Fitt}_A(M') \cdot \operatorname{Fitt}_A(M'') \text{ and } [A : \operatorname{Fitt}_A(M)] = |M|.$$

In particular, if $A = \mathbf{Z}$, then $\operatorname{Fitt}_\mathbf{Z}(M) = |M| \cdot \mathbf{Z}$.

(g) If $A = \mathbf{Z}[G]$, and $G$ is cyclic, then $\operatorname{Fitt}_A(M) \subseteq \operatorname{Fitt}_A(M')$.

Proofs for all the properties above, except (c) and (d), can be found in the Appendix of [9]. (c) is a consequence of (b), and (d) follows easily from the definition of the Fitting ideal.

## 1.5. THE MODULES $X_S$, $U_{S,T}$ AND $\Lambda_{S,T}$

Let $K/k$, $S$, $T$ be as above, and let $r$ be a positive integer. From now on we are going to assume that the set of data $(K/k, S, T, r)$ satisfies the following extended set of hypotheses:

(H) $\begin{cases} S \neq \emptyset, \quad T \neq \emptyset, \quad S \cap T = \emptyset. \\ S \text{ contains all primes which ramify in } K/k. \\ S \text{ contains at least } r \text{ primes which split completely in } K/k. \\ |S| \geq r + 1. \end{cases}$

Hypotheses (H) imply that, for any $\chi \in \widehat{G}$, we have $r_\chi \geq r$ (see (3) above) and therefore $\Theta_{S,T}^{(r)}(0) \in \mathbf{C}[G]$ makes sense. From the definitions, we have $\Theta_{S,T}^{(r)}(0) \in \mathbf{C}[G]_{r,S}$.

Let us choose an $r$-tuple $(v_1, \ldots, v_r)$ of $r$ distinct primes in $S$, which split completely in $K/k$. Let us fix $W = (w_1, \ldots, w_r)$, where $w_i$ is a prime in $K$ lying above $v_i$, for any $1 \leq i \leq r$, and let $w \in S_K$, such that $w \nmid v_i$, for any $1 \leq i \leq r$.

All the exterior powers considered in this paper are over $\mathbf{Z}[G]$, unless otherwise specified. We will be interested in the $\mathbf{Z}[G]$–modules $\overset{r}{\wedge} U_{S,T}$ and $\overset{r}{\wedge} X_S$. Let $\mathbf{x} \overset{\text{def}}{=} (w_1 - w) \wedge \cdots \wedge (w_r - w) \in \overset{r}{\wedge} X_S$.

**Lemma 1.5.1.** *Let $R$ be a discrete valuation ring, or a field, containing $\mathbf{Z}[1/g]$, let $L$ be its field of fractions and let $R[G] = \bigoplus_{\psi \in \widehat{G}(L)} D_\psi$ be the direct sum decomposition described in Remark 1 §1.3. Then*

(1) $RU_{S,T} \xrightarrow{\sim} RX_S \xrightarrow{\sim} \bigoplus_{\psi \in \widehat{G}(L)} D_\psi^{r_\psi}$, *as $R[G]$–modules.*

(2) $\left(R \overset{r}{\wedge} U_{S,T}\right)_{r,S} \xrightarrow{\sim} \left(R \overset{r}{\wedge} X_S\right)_{r,S} \xrightarrow{\sim} R[G]_{r,S}$, *as $R[G]$–modules.*

(3) $\left(R \overset{r}{\wedge} X_S\right)_{r,S} = R[G]_{r,S} \cdot \mathbf{x}$.

**Proof.** (1) is a consequence of Remark 1, §1.3, equality (4) and the fact that $\mathbf{R}U_{S,T}$ and $\mathbf{R}X_S$ are isomorphic as $\mathbf{R}[G]$–modules (see §1.1).

(2) and (3) follow from Remark 1, §1.3 and Lemma 2.6 of [14]. $\square$

For every $\phi_1, \ldots, \phi_r \in \text{Hom}_{\mathbf{Z}[G]}(U_{S,T}, \mathbf{Z}[G])$, one can define a $\mathbf{Q}[G]$–morphism:

$$\mathbf{Q} \overset{r}{\wedge} U_{S,T} \xrightarrow{\phi_1 \wedge \cdots \wedge \phi_r} \mathbf{Q}[G]$$

by letting $\phi_1 \wedge \cdots \wedge \phi_r(u_1 \wedge \cdots \wedge u_r) = \det_{i,j}(\phi_i(u_j))$, for every $u_1 \wedge \cdots \wedge u_r \in \overset{r}{\wedge} U_{S,T}$.

**Definition 1.5.2.** Let $\Lambda_{S,T}$ be the $\mathbf{Z}[G]$–submodule of $\mathbf{Q} \overset{r}{\wedge} U_{S,T}$ defined by:

$$\Lambda_{S,T} = \left\{ \varepsilon \in \left(\mathbf{Q} \overset{r}{\wedge} U_{S,T}\right)_{r,S} \,\middle|\, \begin{array}{l} \phi_1 \wedge \cdots \wedge \phi_r(\varepsilon) \in \mathbf{Z}[G], \\ \forall \phi_1, \ldots, \phi_r \in \text{Hom}_{\mathbf{Z}[G]}(U_{S,T}, \mathbf{Z}[G]) \end{array} \right\}.$$

**Lemma 1.5.3.** *The lattice $\Lambda_{S,T}$ defined above satisfies the equalitites:*

(1) $\mathbf{Z}[1/g]\Lambda_{S,T} = \left(\mathbf{Z}[1/g] \overset{r}{\wedge} U_{S,T}\right)_{r,S}$.

(2) *If $r = 1$, then $\Lambda_{S,T} = (U_{S,T})_{1,S}$.*

(3) *If $r = 0$, then $\Lambda_{S,T} = \mathbf{Z}[G]_{0,S}$.*

**Proof.** See [14], Prop. 1.2. $\square$

## 1.6. THE CONJECTURES

For every component $w_i$ of the $r$–tuple $W$, chosen in §1.5, one can define a $\mathbf{Z}[G]$–morphism:
$$X_S \xrightarrow{w_i^*} \mathbf{Z}[G],$$
by defining it first on $\bigoplus_{w \in S_K} \mathbf{Z}w$ as
$$w_i^*(w) = \sum_{\substack{\sigma \in G \\ w^\sigma = w_i}} \sigma, \quad \forall w \in S_K,$$
and by taking the restricion to $X_S$. One obtains this way a $\mathbf{C}[G]$–morphism
$$\mathbf{C} \overset{r}{\wedge} X_S \xrightarrow{w_1^* \wedge \cdots \wedge w_r^*} \mathbf{C}[G],$$
defined by $w_1^* \wedge \cdots \wedge w_r^*(x_1 \wedge \cdots \wedge x_r) = \det_{i,j}(w_i^*(x_j))$, $\forall x_1 \wedge \cdots \wedge x_r \in \overset{r}{\wedge} X_S$.

**Remark 1.** If $\mathbf{x} = (w_1 - w) \wedge \cdots \wedge (w_r - w)$ is the element of $\overset{r}{\wedge} X_S$, defined in §1.4, then $w_1^* \wedge \cdots \wedge w_r^*(\mathbf{x}) = 1$, as one can easily see.

Let
$$\mathbf{C} \overset{r}{\wedge} U_{S,T} \xrightarrow[\sim]{\lambda_{S,T}^{(r)}} \mathbf{C} \overset{r}{\wedge} X_S$$
be the $\mathbf{C}[G]$–isomorphism defined by $\lambda_{S,T}^{(r)}(u_1 \wedge \cdots \wedge u_r) = \lambda_{S,T}(u_1) \wedge \cdots \wedge \lambda_{S,T}(u_r)$, for every $u_1 \wedge \cdots \wedge u_r \in \overset{r}{\wedge} U_{S,T}$.

**Definition 1.6.1.**
  (1) The regulator map associated to the $r$–tuple $W$ is defined to be the $\mathbf{C}[G]$–morphism:
$$\mathbf{C} \overset{r}{\wedge} U_{S,T} \xrightarrow{R_W} \mathbf{C}[G],$$
  given by $R_W = (w_1^* \wedge \cdots \wedge w_r^*) \circ \lambda_{S,T}^{(r)}$.
  (2) For every $\chi \in \widehat{G}$, we define
$$\mathbf{C} \overset{r}{\wedge} U_{S,T} \xrightarrow{R_{W,\chi}} \mathbf{C}$$
by $R_{W,\chi} = \chi \circ R_W$.

**Remark 2.** Let us observe that:
$$R_W|_{\left(\mathbf{C} \overset{r}{\wedge} U_{S,T}\right)_{r,S}} : \left(\mathbf{C} \overset{r}{\wedge} U_{S,T}\right)_{r,S} \to \mathbf{C}[G]_{r,S}$$
is an isomorphism of $\mathbf{C}[G]$–modules. Indeed, since $\lambda_{S,T}^{(r)}$ is an isomorphism (Lemma 1.5.1), all we have to check is that
$$(w_1^* \wedge \cdots \wedge w_r^*)|_{\left(\mathbf{C} \overset{r}{\wedge} X_S\right)_{r,S}} : \left(\mathbf{C} \overset{r}{\wedge} X_S\right)_{r,S} \to \mathbf{C}[G]_{r,S}$$
is an isomorphism. Acording to Lemma 1.5.1 (3), $\left(\mathbf{C} \overset{r}{\wedge} X_S\right)_{r,S} = \mathbf{C}[G]_{r,S} \cdot \mathbf{x}$, and therefore our statement follows from Remark 1 above.

**Conjecture A ("over Q")(Stark, Tate).** *If the set of data $(K/k, S, T, r)$ satisfies hypotheses* (H), *then there exists a unique $\varepsilon_{S,T} \in \left(\mathbf{Q} \overset{r}{\wedge} U_{S,T}\right)_{r,S}$, such that $R_W(\varepsilon_{S,T}) = \Theta_{S,T}^{(r)}(0)$.*

**Remark 3.** As Rubin shows (see [14], Prop. 2.3), Conjecture A above is equivalent to Stark's original conjecture "over $\mathbf{Q}$", for the $L$–functions associated to characters $\chi \in \widehat{G}(\mathbf{C}, r)$. As Tate points out (see [17], Chpt.V), Stark's Conjecture "over $\mathbf{Q}$" is always true in the function field case. Therefore Conjecture A is a theorem in the setting we are considering in this paper. We will refer to Conjecture A as $\mathbf{Q} \cdot \mathrm{St}\,(K/k, S, T, r)$ in the sequel.

The following is the refinement Rubin [14] proposes for Conjecture A:

**Conjecture B ("over Z") (Rubin).** *If the set of data $(K/k, S, T, r)$ satisfies hypotheses* (H), *then there exists a unique $\varepsilon_{S,T} \in \Lambda_{S,T}$, such that $R_W(\varepsilon_{S,T}) = \Theta_{S,T}^{(r)}(0)$.*

**Remark 4.** If $r = 0$, this statement is equivalent to $\Theta_{S,T}(0) \in \mathbf{Z}[G]_{0,S}$ (see Lemma 1.5.3(3)). This was proved by Deligne in the function field case ( see [17], Chpt.V).

If $r = 1$, as Rubin points out (see [14], Prop. 2.5), Conjecture B above, for a fixed $S$ and for all $T$, is equivalent to Stark's own refined conjecture for all $L$–functions associated to characters $\chi \in \widehat{G}(\mathbf{C}, r)$. In the function field case, this was proved independently and with very different methods by Deligne (see [17], Chpt.V) and Hayes [6]. (Also see [13 §3] for more details on the equivalence between Rubin's Conjecture for $r = 1$ and the classical Brumer–Stark Conjecture.)

The main goal of this paper is the study of Conjecture B above (which will be refered to as $\mathrm{St}\,(K/k, S, T, r)$ ), for any value of $r$.

**Remark 5.** Remark 2 above implies that the uniqueness of $\varepsilon_{S,T}$ in both Conjecture A and B is automatic, once one proves its existence in the appropriate vector space or lattice.

1.7. GEOMETRIC BACKGROUND (see [17])

**1.7.1 The corresponding schemes.** Let $X_1 \longrightarrow Y_1$ be a finite morphism of projective, irreducible, smooth curves, defined over $\mathbf{F}_q$, corresponding to the inclusion $k \hookrightarrow K$. Let $\mathbf{F}$ be an algebraic closure of $\mathbf{F}_q$ and let $X = X_1 \underset{\mathrm{Spec}(\mathbf{F}_q)}{\times} \mathrm{Spec}(\mathbf{F})$, and $Y = Y_1 \underset{\mathrm{Spec}(\mathbf{F}_q)}{\times} \mathrm{Spec}(\mathbf{F})$ be the smooth, projective curves associated to $X_1$ and $Y_1$ respectively, by extending scalars from $\mathbf{F}_q$ to $\mathbf{F}$. Since $\mathbf{F}_q$ is algebraically closed in $k$, $k \underset{\mathbf{F}_q}{\otimes} \mathbf{F}$ is a field and, correspondingly, $Y$ is an irreducible, smooth, projective curve, defined over $\mathbf{F}$. Since the algebraic closure of $\mathbf{F}_q$ in $K$ is $\mathbf{F}_{q^\nu}$, $K \underset{\mathbf{F}_q}{\otimes} \mathbf{F} \overset{\sim}{\longrightarrow} \underset{0 \leq i \leq \nu - 1}{\bigoplus} K^{(i)}$, with $K^{(i)}$ mutually isomorphic fields. Correspondingly, $X = \coprod_{0 \leq i \leq \nu - 1} X^{(i)}$, where the $X^{(i)}$'s are irreducible, smooth, projective curves, defined over $\mathbf{F}$, whose fields of rational functions are the $K^{(i)}$'s respectively. The $X^{(i)}$'s are mutually isomorphic over $\mathbf{F}$, and therefore their genera $g_{X^{(i)}}$ are equal. We will denote the common value of these numbers by $g_X$. We will also denote by $[X^{(i)}]$ the generic point of the scheme $X^{(i)}$, for every $0 \leq i \leq \nu - 1$.

Every $\sigma \in G$ gives isomorphisms

$$\sigma^{-1} : K \xrightarrow{\sim} K, \quad \sigma^{-1} \otimes 1_\mathbf{F} : K \otimes \mathbf{F} \xrightarrow{\sim} K \otimes \mathbf{F},$$

which induce the isomorphisms

$$\sigma_{X_1} : X_1 \xrightarrow{\sim} X_1, \quad \sigma_X : X \xrightarrow{\sim} X$$

of $\mathbf{F}_q$–schemes and $\mathbf{F}$–schemes respectively. This way, $G$ acts on the sets of points of $X$ and it permutes the generic points $[X^{(i)}]$ transitively. We emphasize the fact that $\sigma_{X_1}$ and $\sigma_X$ are associated to $\sigma^{-1} \in G$. This way, if $w$ is a prime in $K$, and $[w]$ is the associated closed point on the scheme $X_1$, then $\sigma_{X_1}([w]) = [w^\sigma]$.

**1.7.2 The geometric Frobenius.** Let $F_{X_1} : X_1 \longrightarrow X_1$ be the geometric Frobenius endomorphism of $X_1$, relative to $\mathbf{F}_q$. If $U$ is an affine open subset of $X_1$, then $F_{X_1}|_U$ corresponds to the $q$–power map on the $\mathbf{F}_q$–algebra $\Gamma(U, O_{X_1})$, of $U$-sections of the structural sheaf $O_{X_1}$.

**Definition 1.7.2.** Let

$$F = F_{X_1} \underset{\mathrm{Spec}(\mathbf{F}_q)}{\times} 1_{\mathrm{Spec}(\mathbf{F})} : X \longrightarrow X$$

be the $\mathrm{Spec}(\mathbf{F})$–scheme morphism obtained from $F_{X_1}$ by extending scalars to $\mathbf{F}$. $F$ is called the geometric Frobenius endomorphism of $X$ relative to $\mathbf{F}_q$.

**Remark 1.** $F$ permutes transitively the generic points $[X^{(i)}]$ of $X$, and that $F^\nu$ fixes each of these. We can therefore suppose from now on that $[X^{(i)}] = F^i([X^{(0)}])$, for every $0 \leq i \leq \nu - 1$.

**Remark 2.** The simple fact that $\sigma(x^q) = \sigma(x)^q, \forall x \in K, \forall \sigma \in G$, implies that $F_{X_1} \circ \sigma_{X_1} = \sigma_{X_1} \circ F_{X_1}$ and therefore, by extending scalars, $F \circ \sigma_X = \sigma_X \circ F, \forall \sigma \in G$.

**1.7.3 The $\ell$–adic homology groups of $X$.** We first give explicit descriptions of the homology groups $\mathrm{H}_i(X, \mathbf{Z})$ of $X$, with coefficients in $\mathbf{Z}$, and of the actions of $G$ and of the geometric Frobenius map $F$ on each of these. For any $i \geq 0$ and any $\sigma \in G$, we will denote by $\sigma_{*,i}$ and $F_{*,i}$ the maps induced on $\mathrm{H}_i(X, \mathbf{Z})$ by $\sigma$ and $F$ respectively.

- $\mathrm{H}_i(X, \mathbf{Z}) = 0$ for $i \geq 3$.

- $H_0(X, \mathbf{Z}) =$ the free abelian group generated by $\left[X^{(i)}\right], 0 \leq i \leq \nu - 1$.

The maps $\sigma_{*,0}$ and $F_{*,0}$ act on $H_0(X, \mathbf{Z})$ by simply permuting the $[X^{(i)}]$'s as described in §1.7.2. Remark 1 in §1.7.2 also implies that $F_{*,0}^\nu = 1_{H_0(X, \mathbf{Z})}$.

- $H_1(X, \mathbf{Z}) = \mathrm{Pic}^0(X)$, where $\mathrm{Pic}^0(X)$ is the Jacobian variety associated to $X$. We have $\mathrm{Pic}^0(X) = \coprod_i \mathrm{Pic}^0(X^{(i)})$, and each $\mathrm{Pic}^0(X^{(i)})$ is an irreducible abelian variety, of dimension $g_X$, whose underlying group is isomorphic to the quotient of the group of divisors of degree zero on $X^{(i)}$, by the subgroup of principal divisors:

$$\mathrm{Pic}^0(X^{(i)}) = \frac{\mathrm{Div}^0\left(X^{(i)}\right)}{\{\mathrm{div}(f) : f \in K^{(i)\times}\}}, \quad \mathrm{Pic}^0(X) = \frac{\mathrm{Div}^0(X)}{\{\mathrm{div}(f) : f \in \bigoplus_i K^{(i)\times}\}}.$$

Here $\mathrm{Div}^0(X)$ is the group of divizors of *multidegree* zero on $X$, i.e. of degree zero if restricted to each component $X^{(i)}$.

The maps $\sigma_{*,1}$ and $F_{*,1}$ are naturally induced by the actions of $\sigma$ and $F$ on the set of closed points of the scheme $X$.

- $H_2(X, \mathbf{Z}) = H_0(X, \mathbf{Z}) \otimes \mathbf{F}^\times$, with $\sigma_{*,2} = \sigma_{*,0} \otimes 1_{\mathbf{F}^\times}$, and $F_{*,2} = F_{*,0} \otimes \sigma_q$, where $\sigma_q : \mathbf{F}^\times \longrightarrow \mathbf{F}^\times$ is the $q$–power map $\sigma_q(x) = x^q$.

If $M$ is an abelian group and $\ell$ is a prime number, we give the following ad–hoc definition of the $\ell$–adic Tate module of $M$:

$$T_\ell(M) = \begin{cases} M \otimes \mathbf{Z}_\ell & \text{, if } M \text{ has no } \mathbf{Z}\text{–torsion} \\ \varprojlim_n M[\ell^n] & \text{, otherwise,} \end{cases}$$

where $M[\ell^n]$ is the group of $\ell^n$–torsion points of $M$, and the projective limit is taken with respect to the multiplication by $\ell$ maps $M[\ell^{n+1}] \xrightarrow{\ell} M[\ell^n]$.

**Definition 1.7.3.1.** Let $\ell$ be any prime number. The $\ell$–adic homology groups $H_i(X, \mathbf{Z}_\ell)$ of $X$ are defined by:

$$H_i(X, \mathbf{Z}_\ell) = T_\ell(H_i(X, \mathbf{Z})), \quad \forall i \geq 0.$$

The $\mathbf{Z}$–linear actions of $G$ and $F$ on $H_i(X, \mathbf{Z})$ induce in a natural way $\mathbf{Z}_\ell$–linear actions on $H_i(X, \mathbf{Z}_\ell)$. According to the definition above, we have:

- $H_i(X, \mathbf{Z}_\ell) = 0, \forall i \geq 3$.
- $H_0(X, \mathbf{Z}_\ell) = H_0(X, \mathbf{Z}) \otimes \mathbf{Z}_\ell$ is the free $\mathbf{Z}_\ell$–module generated by $[X^{(i)}]$, for all $0 \leq i \leq \nu - 1$.
- $H_1(X, \mathbf{Z}_\ell) = T_\ell\left(\mathrm{Pic}^0(X)\right) = \bigoplus_i T_\ell\left(\mathrm{Pic}^0(X^{(i)})\right) \xrightarrow{\sim} \mathbf{Z}_\ell^{\nu n}$, where $n = 2g_X$, if $\ell \neq p$, and $0 \leq n < 2g_X$, if $\ell = p$.
- $H_2(X, \mathbf{Z}_\ell) \xrightarrow{\sim} H_0(X, \mathbf{Z}_\ell) \otimes T_\ell(\mathbf{F}^\times)$. Since $T_\ell(\mathbf{F}^\times) = \varprojlim \mu_{\ell^n}$, where $\mu_{\ell^n} = \{\zeta \in \mathbf{F}^\times : \zeta^{\ell^n} = 1\}$, we have $T_\ell(\mathbf{F}^\times) \xrightarrow{\sim} \mathbf{Z}_\ell$ if $\ell \neq p$, and $T_p(\mathbf{F}^\times) = \{1\}$. We therefore have $H_2(X, \mathbf{Z}_\ell) \xrightarrow{\sim} H_0(X, \mathbf{Z}_\ell)$, if $\ell \neq p$, and $H_2(X, \mathbf{Z}_p) = 0$. It is worth noticing that the isomorphism $H_2(X, \mathbf{Z}_\ell) \xrightarrow{\sim} H_0(X, \mathbf{Z}_\ell)$ is $G$–equivariant, but not $F_*$–equivariant. In fact, the action of $F_{*,2}$ is taken into the action of $q \cdot F_{*,0}$ by this isomorphism.

If $\ell$ is a prime number, and $R$ is any ring containing $\mathbf{Z}_\ell$, then, by definition, $\mathrm{H}_i(X, R) = \mathrm{H}_i(X, \mathbf{Z}_\ell) \underset{\mathbf{Z}_\ell}{\bigotimes} R$. In what follows, $\sigma_{*,i}$ and $F_{*,i}$ will also denote the $R$–linear actions induced by $\sigma$ and $F$ on $\mathrm{H}_i(X, R)$, for any $i \geq 0$, and any $R$ as above.

**Remark 1.** The $\mathrm{H}_i(X, \mathbf{Z}_\ell)$'s as defined above are in fact the functorial duals of the étale cohomology groups $\mathrm{H}^i_{\text{ét}}(X, \mathbf{Z}_\ell)$ of $X$ (see Appendix).

**1.7.4 Homological interpretation of $\Theta_{S,T}(s)$.** Let $R$ be a commutative ring, and $V$ a finitely generated projective $R[G]$–module. Let $W$ be a finitely generated $R[G]$–module so that $V \oplus W \xrightarrow{\sim} R[G]^n$, for some integer $n \geq 0$. Let $f \in \mathrm{End}_{R[G]}(V)$. If $u$ is a variable, one can define the polynomial in $u$, with coefficients in $R[G]$:

$$\mathrm{det}_{R[G]}(1 - f \cdot u \mid V) \stackrel{\text{def}}{=} \mathrm{det}_{R[G]}(1 - (f \oplus 0_W) \cdot u \mid V \oplus W),$$

where the determinant on the right is taken with respect to a basis of the free $R[G]$–module $V \oplus W$. Schanuel's Lemma easily implies that this definition does not depend on $W$.

In particular, if $R$ and $L$ are as in Remark 1, §1.3, then one can consider the decompositions $R[G] = \underset{\psi \in \widehat{G}(L)}{\bigoplus} D_\psi$, and $V = \underset{\psi \in \widehat{G}(L)}{\bigoplus} V^\psi$, with $V^\psi$ free $D_\psi$–modules of finite rank. It is easy to see that, in this case,

$$\mathrm{det}_{R[G]}(1 - f \cdot u \mid V) = \sum_{\psi \in \widehat{G}(L)} \mathrm{det}_{D_\psi}(1 - f \cdot u \mid V^\psi),$$

in $R[G][u] = \underset{\psi \in \widehat{G}(L)}{\bigoplus} D_\psi[u]$.

**Remark 1.** If $R'$ is an $R$–algebra, $V$ is a projective $R$–module, and $f \in \mathrm{End}_{R[G]}(V)$, then $V' = V \otimes_R R'$ is a projective $R'[G]$–module, $f' = f \otimes 1_{R'} \in \mathrm{End}_{R'[G]}(V')$, and $\mathrm{det}_{R'[G]}(1 - f' \cdot u \mid V') = \mathrm{det}_{R[G]}(1 - f \cdot u \mid V)$.

For every prime number $\ell$, and every $i = 0, 1, 2$, let

$$P_{i,\ell}(u) \stackrel{\text{def}}{=} \mathrm{det}_{\mathbf{Q}_\ell[G]}(1 - F_{*,i} \cdot u \mid \mathrm{H}_i(X, \mathbf{Q}_\ell)).$$

The following theorem will play an important role in our future arguments.

**Theorem 1.7.4.1.** *For any $i = 0, 1, 2$:*
 (1) $P_{i,p}(u) \in \mathbf{Z}_p[1/g][G][u]$, *and if $\ell \neq p$, then $P_{i,\ell} \in \mathbf{Z}[1/g][G][u]$.*
 (2) *If $\ell \neq p$, then $P_{i,\ell}(u)$ does not depend on $\ell$.*
   *Let $P_i(u) \stackrel{\text{def}}{=} P_{i,\ell}(u)$, for any $\ell \neq p$. Then:*
 (3) *There exist polynomials $Q_i(u) \in \mathbf{Z}_p[1/g][G][u]$, such that*

$$P_i(u) = P_{i,p}(u) \cdot Q_i(u).$$

 (4) *If $K/k$ is a constant field extension, then $Q_i(u) \in \mathbf{Z}_p[G][u]$.*

**Proof.** See Appendix. □

The link between the polynomials $P_i(u)$ defined above and the Stickelberger function $\Theta(s)$ is given by the following:

**Theorem 1.7.4.2.** *There is an equality of complex, meromorphic functions*

$$\Theta(s) = \prod_{0 \leq i \leq 2} P_i(q^{-s})^{(-1)^{i+1}}.$$

**Proof.** See [17], Chpt.V. □

The natural question which arises is to give a similar homological interpretation for $\Theta_{S,T}(s)$. This was done by Deligne (see [17], Chpt.V). We are going to briefly describe his constructions below, as they will be needed in §§3 and 4.

Let $S_X$ be the finite set of closed points on $X$, lying above points corresponding to primes in $S$, and let $\text{Div}(S_X)$ be the $\mathbf{Z}[G]$–module of divisors on $X$, supported on $S_X$. There is an exact sequence of $\mathbf{Z}[G]$–modules:

$$0 \longrightarrow \ker \alpha_S \longrightarrow \text{Div}(S_X) \xrightarrow{\alpha_S} \text{H}_0(X, \mathbf{Z}) \longrightarrow 0, \tag{5}$$

where $\alpha_S$ is the multidegree map given by:

$$\alpha_S \left( \sum_{j=1}^m n_j P_j \right) = \sum_{0 \leq i \leq \nu-1} \left( \sum_{P_j \in X^{(i)}} n_j \right) \left[ X^{(i)} \right],$$

for every $\sum_{j=1}^m n_j P_j \in \text{Div}(S_X)$.

**Remark 2.** Obviously $\ker \alpha_S$ is a free $\mathbf{Z}$–module of finite rank and, by definition, $\text{T}_\ell(\ker \alpha_S) = \ker \alpha_S \otimes \mathbf{Z}_\ell$, for any prime number $\ell$. $\alpha_S$ is $F_{*,0}$–invariant and therefore $F$ acts on both $\ker \alpha_S$ and $\text{T}_\ell(\ker \alpha_S)$. We denote the induced endomorphisms by $F_{*,0}$.

Let $T_X$ be the finite set of closed points on $X$, lying above points associated to primes in $T$, and let $\text{H}_0(T_X, \mathbf{Z})$ be the free abelian group generated by $T_X$. There is an exact sequence of $\mathbf{Z}[G]$–modules

$$0 \longrightarrow \text{H}_0(X, \mathbf{Z}) \otimes \mathbf{F}^\times \xrightarrow{\beta_T} \text{H}_0(T_X, \mathbf{Z}) \otimes \mathbf{F}^\times \longrightarrow \text{coker}(\beta_T) \longrightarrow 0, \tag{6}$$

where

$$\beta_T \left( \sum_{0 \leq i \leq \nu-1} \left[ X^{(i)} \right] \otimes f_i \right) = \sum_{0 \leq i \leq \nu-1} \left( \sum_{P \in T_X \cap X^{(i)}} P \right) \otimes f_i,$$

for any $\sum_{0 \leq i \leq \nu-1} \left[ X^{(i)} \right] \otimes f_i \in \text{H}_0(X, \mathbf{Z}) \otimes \mathbf{F}^\times$.

**Remark 3.** From the definition of $\beta_T$, one can easily see that, as a group, $\text{coker}(\beta_T)$ is a torus, i.e. $\text{coker}(\beta_T) \xrightarrow{\sim} \mathbf{F}^{\times m}$, for some positive integer $m$. Therefore we have $\text{T}_\ell(\text{coker}(\beta_T)) \xrightarrow{\sim} \mathbf{Z}_\ell^m$, if $\ell \neq p$, and $\text{T}_p(\text{coker}(\beta_T)) = \{1\}$.

Since $F$ induces an isomorphism $F_{*,2} = F_{*,0} \otimes \sigma_q$ on both $\mathrm{H}_0(X, \mathbf{Z}) \otimes \mathbf{F}^\times$ and $\mathrm{H}_0(T_X, \mathbf{Z}) \otimes \mathbf{F}^\times$, and $\beta_T$ is $F_{*,2}$–invariant, $F$ induces a $\mathbf{Z}$–linear isomorphism on $\mathrm{coker}(\beta_T)$, and a $\mathbf{Z}_\ell$–linear isomorphism on $\mathrm{T}_\ell(\mathrm{coker}(\beta_T))$. We will denote both these isomorphisms by $F_{*,2}$.

**Theorem 1.7.4.3.** *If $S$ and $T$ are two finite, nonempty, disjoint sets of primes in $k$, $S$ containing all primes which ramify in $K/k$, and $u = q^{-s}$, then*

(1) $\Theta_{S,T}(s) \in \mathbf{Z}[G][u]$.
(2) *For every prime number $\ell$, such that $\ell \neq p$, we have*

$$\Theta_{S,T}(s) = \mathrm{det}_{\mathbf{Q}_\ell[G]}\left(1 - F_{*,0} \cdot u \mid \mathrm{T}_\ell(\ker \alpha_S) \otimes \mathbf{Q}_\ell\right) \cdot$$
$$\cdot \mathrm{det}_{\mathbf{Q}_\ell[G]}\left(1 - F_{*,1} \cdot u \mid \mathrm{H}_1(X, \mathbf{Q}_\ell)\right) \cdot$$
$$\cdot \mathrm{det}_{\mathbf{Q}_\ell[G]}\left(1 - F_{*,2} \cdot u \mid \mathrm{T}_\ell(\mathrm{coker}(\beta_T)) \otimes \mathbf{Q}_\ell\right).$$

**Proof.** See [17], Chpt.V. □

**Remark 4.** The theorem above implies in particular that, if the set of data $(K/k, S, T, 0)$ satisfies hypotheses (H), then $\Theta_{S,T}(0) \in \mathbf{Z}[G]$, and therefore (by definition) $\Theta_{S,T}(0) \in \mathbf{Z}[G]_{0,S}$. This is the proof of $\mathrm{St}(K/k, S, T, r)$, for $r = 0$.

**Lemma 1.7.4.4.** *Let $\ell$ be a prime number, such that $\gcd(\ell, g) = 1$. Then*

(1) *If $\ell \neq p$,*

$$\Theta_{S,T}(s) = \mathrm{det}_{\mathbf{Z}_\ell[G]}\left(1 - F_{*,0} \cdot u \mid \mathrm{T}_\ell(\ker \alpha_S)\right) \cdot \mathrm{det}_{\mathbf{Z}_\ell[G]}\left(1 - F_{*,1} \cdot u \mid \mathrm{H}_1(X, \mathbf{Z}_\ell)\right)$$
$$\cdot \mathrm{det}_{\mathbf{Z}_\ell[G]}\left(1 - F_{*,2} \cdot u \mid \mathrm{T}_\ell(\mathrm{coker}(\beta_T))\right).$$

(2) *If $\gcd(p, g) = 1$, there exists a polynomial $Q(u) \in \mathbf{Z}_p[G][u]$ such that*

$$\Theta_{S,T}(s) = Q(u) \cdot \mathrm{det}_{\mathbf{Z}_p[G]}\left(1 - F_{*,1} \cdot u \mid \mathrm{H}_1(X, \mathbf{Z}_p)\right) \cdot$$
$$\cdot \mathrm{det}_{\mathbf{Z}_p[G]}\left(1 - F_{*,0} \cdot u \mid \mathrm{T}_p(\ker \alpha_S)\right).$$

**Proof.** Let us first notice that, for $\ell$ satisfying $\gcd(\ell, g) = 1$, $\mathrm{H}_1(X, \mathbf{Z}_\ell)$, $\mathrm{T}_\ell(\ker \alpha_S)$ and $\mathrm{T}_\ell(\mathrm{coker}(\beta_T))$ are $\mathbf{Z}_\ell$–free and therefore $\mathbf{Z}_\ell[G]$–projective modules (see Remark 1, §1.3). This shows that the determinants on the right–hand side of (1) make sense. Statement (1) in the lemma is now a direct consequence of Remark 1 above and Theorem 1.7.4.3(2).

The fact that $\mathrm{T}_\ell(\ker \alpha_S) = (\ker \alpha_S \otimes \mathbf{Z}[1/g]) \bigotimes_{\mathbf{Z}[1/g]} \mathbf{Z}_\ell$ and that $\ker \alpha_S \otimes \mathbf{Z}[1/g]$ is $\mathbf{Z}[G]$–projective imply that

$$\mathrm{det}_{\mathbf{Z}[1/g][G]}(1 - F_{*,0} \cdot u \mid \ker \alpha_S \otimes \mathbf{Z}[1/g]) = \mathrm{det}_{\mathbf{Z}_\ell[G]}\left(1 - F_{*,0} \cdot u \mid \mathrm{T}_\ell(\ker \alpha_S)\right), \quad (7)$$

for all prime numbers $\ell$, such that $\gcd(\ell, g) = 1$, in particular for $\ell = p$. Equality (7) shows that the polynomial on the right is independent of $\ell$, as long as $\gcd(\ell, g) = 1$.

As Deligne shows (see [17], Chpt.V), the exact sequence (6) implies that, if $\ell \neq p$, $\gcd(\ell, g) = 1$, we have an equality

$$\det\nolimits_{\mathbf{Z}_\ell[G]} \left(1 - F_{*,2} \cdot u \mid \mathrm{T}_\ell(\mathrm{coker}(\beta_T))\right) = \frac{\prod\limits_{v \in T} (1 - \sigma_v^{-1}(qu)^{d_v})}{\det\nolimits_{\mathbf{Z}_\ell[G]} \left(1 - qF_{*,0} \cdot u \mid \mathrm{H}_0\left(X, \mathbf{Z}_\ell\right)\right)}. \quad (8)$$

On the other hand, $\mathrm{H}_0\left(X, \mathbf{Z}_\ell\right) = \mathrm{H}_0\left(X, \mathbf{Z}\left[1/g\right]\right) \otimes \mathbf{Z}_\ell$, and $\mathrm{H}_0\left(X, \mathbf{Z}\left[1/g\right]\right)$ is a projective $\mathbf{Z}\left[1/g\right][G]$–module. Remark 1 therefore implies that

$$\det\nolimits_{\mathbf{Z}_\ell[G]} \left(1 - qF_{*,0} \cdot u \mid \mathrm{H}_0\left(X, \mathbf{Z}_\ell\right)\right) = \det\nolimits_{\mathbf{Z}[1/g][G]} \left(1 - qF_{*,0} \cdot u \mid \mathrm{H}_0\left(X, \mathbf{Z}\left[1/g\right]\right)\right).$$

This implies that the right hand–side of equality (8) is an element in the power series ring $\mathbf{Z}\left[1/g\right][G][[u]]$. On the other hand, the left–hand side of the same equality belongs to $\mathbf{Z}_\ell[G][u]$. We therefore have

$$\det\nolimits_{\mathbf{Z}_\ell[G]} \left(1 - F_{*,2} \cdot u \mid \mathrm{T}_\ell(\mathrm{coker}(\beta_T))\right) \in \mathbf{Z}\left[1/g\right][G][u],$$

for all $\ell$, such that $\gcd(\ell, g) = 1$.

Let $Q(u) \stackrel{def}{=} Q_1(u) \cdot \det\nolimits_{\mathbf{Z}_\ell[G]} \left(1 - F_{*,2} \cdot u \mid \mathrm{T}_\ell(\mathrm{coker}(\beta_T))\right)$, where $Q_1(u)$ is the polynomial defined in Theorem 1.7.4.1 (3), and $\ell$ is a prime number, such that $\ell \neq p$ and $\gcd(\ell, g) = 1$. According to the arguments above, $Q(u) \in \mathbf{Z}_p[G][u]$. Part (1) of the lemma, in addition to Theorem 1.7.4.1(3) and equality (7) imply that

$$\Theta_{S,T}(s) = Q(u) \cdot \det\nolimits_{\mathbf{Z}_p[G]} \left(1 - F_{*,1} \cdot u \mid \mathrm{H}_1\left(X, \mathbf{Z}_p\right)\right) \cdot$$
$$\cdot \det\nolimits_{\mathbf{Z}_p[G]} \left(1 - F_{*,0} \cdot u \mid \mathrm{T}_p(\ker \alpha_S)\right),$$

which concludes the proof of (2). $\square$

## 2. Dynamics (Changing $S$ and $r$)

Let us suppose that the set of data $(K/k, S, T, r)$ satisfies hypotheses (H). Let $r' > r$, and let $v_{r+1}, \ldots, v_{r'}$ be $r' - r$ distinct primes in $k$, not belonging to $S$ or $T$, and splitting completely in $K/k$. If $S' = S \cup \{v_{r+1}, \ldots, v_{r'}\}$, then the set of data $(K/k, S', T, r')$ satisfies hypotheses (H) as well. Since $T$ remains fixed throughout this section, we will denote by $\varepsilon_S$ and $\varepsilon_{S'}$ the elements $\varepsilon_{S,T}$ and $\varepsilon_{S',T}$ whose existence is predicted by $\mathbf{Q} \cdot \mathrm{St}\,(K/k, S, T, r)$ and $\mathbf{Q} \cdot \mathrm{St}\,(K/k, S', T, r')$ respectively. We will also make the notations $A_S = A_{S,T}$, $U_S = U_{S,T}$ etc. $A_{S,S'}$ will denote the subgroup of $A_S$ generated by primes in $K$, lying above $v_i$, for all $r + 1 \leq i \leq r'$. All the Fitting ideals and exterior powers involved from now on in this paper are considered over $\mathbf{Z}[G]$, unless otherwise specified, and therefore we will suppress the group ring $\mathbf{Z}[G]$ from the notation.

In this section we will list some relations between statements $\mathrm{St}\,(K/k, S, T, r)$ and $\mathrm{St}\,(K/k, S', T, r')$, which will be needed in §§3 and 4. For most proofs we refer the reader to [14], §5.

There is an exact sequence of $\mathbf{Z}[G]$–modules:

$$0 \longrightarrow A_{S,S'} \longrightarrow A_S \longrightarrow A_{S'} \longrightarrow 0. \quad (9)$$

**Proposition 2.1.** *If $\varepsilon_S$ and $\varepsilon_{S'}$ are the elements defined above, then*
  (1) $\varepsilon_S \in \mathbf{Z}[1/g]\operatorname{Fitt}(A_S)\Lambda_S \iff \varepsilon_{S'} \in \mathbf{Z}[1/g]\operatorname{Fitt}(A_{S'})\Lambda_{S'}$.
  (2) *If $G$ is cyclic, then $\varepsilon_S \in \operatorname{Fitt}(A_S)\Lambda_S \Longrightarrow \varepsilon_{S'} \in \Lambda_{S'}$.*

**Proof.** (see also Theorem 5.3 in [14]) According to [14], §5, there exist $\mathbf{Q}[G]$–morphisms:
$$\mathbf{Q}\overset{r'}{\wedge} U_{S'} \overset{\Phi}{\longrightarrow} \mathbf{Q}\overset{r}{\wedge} U_S, \qquad \mathbf{Q}\overset{r'-r}{\wedge} U_{S'} \overset{\Phi'}{\longrightarrow} \mathbf{Q}[G],$$
satisfying the following properties:
  (a) $\Phi$ is injective if restricted to $\mathbf{Q}\Lambda_{S'}$.
  (b) $\Phi(\mathbf{Z}[1/g]\Lambda_{S'}) = \mathbf{Z}[1/g]\operatorname{Fitt}(A_{S,S'})\Lambda_S$.
  (c) $\operatorname{Fitt}(A_{S,S'})\Lambda_S \subseteq \Phi(\Lambda_{S'})$.
  (d) $\Phi(\varepsilon_{S'}) = \varepsilon_S$.
  (e) $\Phi'\left(\overset{r'-r}{\wedge} U_{S'}\right) = \operatorname{Fitt}(A_{S,S'})$.
  (f) $\Phi(u_1 \wedge u_2) = \Phi'(u_1) \wedge u_2$, for every $u_1 \in \mathbf{Q}\overset{r'-r}{\wedge} U_{S'}$ and $u_2 \in \mathbf{Q}\overset{r}{\wedge} U_S$.

(1) The exact sequence (9) and the fact that $\mathbf{Z}[1/g][G]$ is a finite, direct sum of Dedekind domains (see §1.3) imply that there is an equality of $\mathbf{Z}[1/g][G]$–ideals (see §1.4(f)):
$$\mathbf{Z}[1/g]\operatorname{Fitt}(A_S) = (\mathbf{Z}[1/g]\operatorname{Fitt}(A_{S,S'})) \cdot (\mathbf{Z}[1/g]\operatorname{Fitt}(A_{S'})). \qquad (10)$$

Let us suppose that $\varepsilon_S \in \mathbf{Z}[1/g]\operatorname{Fitt}(A_S)\Lambda_S$. (10) implies that $\varepsilon_S$ can be written as a finite sum $\varepsilon_S = \sum_i \gamma_i \delta_i \lambda_i$, with $\gamma_i \in \mathbf{Z}[1/g]\operatorname{Fitt}(A_{S'})$, $\delta_i \in \operatorname{Fitt}(A_{S,S'})$ and $\lambda_i \in \Lambda_S$. According to (c) above, this implies that $\varepsilon_S \in \Phi(\mathbf{Z}[1/g]\operatorname{Fitt}(A_{S'})\Lambda_{S'})$. Thus (a) and (d) imply that $\varepsilon_{S'} \in \mathbf{Z}[1/g]\operatorname{Fitt}(A_{S'})\Lambda_{S'}$.

Let us suppose now that $\varepsilon_{S'} \in \mathbf{Z}[1/g]\operatorname{Fitt}(A_{S'})\Lambda_{S'}$. This implies that $\varepsilon_S = \Phi(\varepsilon_{S'}) \in \operatorname{Fitt}(A_{S'})\Phi(\mathbf{Z}[1/g]\Lambda_{S'})$. On the other hand, according to (b) and (10) above, we have the equalities:
$$\operatorname{Fitt}(A_{S'})\Phi(\mathbf{Z}[1/g]\Lambda_{S'}) = \mathbf{Z}[1/g]\operatorname{Fitt}(A_{S,S'})\operatorname{Fitt}(A_{S'})\Lambda_S =$$
$$= \mathbf{Z}[1/g]\operatorname{Fitt}(A_S)\Lambda_S.$$

This concludes the proof of (1).
  (2) If $G$ is cyclic, §1.4 (g) and the exact sequence (9) imply that $\operatorname{Fitt}(A_S) \subseteq \operatorname{Fitt}(A_{S,S'})$, and therefore (a) and (c) imply (2). $\square$

**Corollary 2.2.** *Let us suppose that, in the context above, $r = 0$. Then:*
  (1) $\Theta_{S,T}(0) \in \mathbf{Z}[1/g]\operatorname{Fitt}(A_S) \iff \varepsilon_{S'} \in \mathbf{Z}[1/g]\operatorname{Fitt}(A_{S'})\Lambda_{S'}$.
  (2) *If $G$ is cyclic, $\Theta_{S,T}(0) \in \mathbf{Z}[G]_{0,S}\operatorname{Fitt}(A_S) \Longrightarrow \varepsilon_{S'} \in \mathbf{Z}[G]_{r',S'} \cdot \overset{r'}{\wedge} U_{S'}$.*

**Proof.** (1) In order to prove (1), one has to rely on the fact that, in the case $r = 0$, $\Lambda_S = \mathbf{Z}[G]_{0,S}$ (see Lemma 1.5.3(3)), and that

$$\mathbf{Z}[1/g]\operatorname{Fitt}(A_S) \cdot \mathbf{Z}[1/g][G]_{0,S} = \mathbf{Z}[1/g]\operatorname{Fitt}(A_S) \cap \mathbf{Z}[1/g][G]_{0,S}.$$

This last equality follows from the fact that $\mathbf{Z}[1/g][G] = \bigoplus_{\psi \in \widehat{G}(\mathbf{Q})} D_\psi$ and that, with respect to this decomposition, $\mathbf{Z}[1/g][G]_{0,S} = \bigoplus_{\psi \in \widehat{G}(\mathbf{Q},0)} D_\psi$ (see §1.3).

Proposition 2.1 (1), in addition to the fact that $\Theta_{S,T}(0) \in \mathbf{C}[G]_{0,S}$, implies therefore that

$$\Theta_{S,T}(0) \in \mathbf{Z}[1/g]\operatorname{Fitt}(A_S) \Longleftrightarrow \Theta_{S,T}(0) \in \mathbf{Z}[1/g]\operatorname{Fitt}(A_S) \cdot \mathbf{Z}[1/g][G]_{0,S}$$
$$\Longleftrightarrow \varepsilon_{S'} \in \mathbf{Z}[1/g]\operatorname{Fitt}(A_{S'})\Lambda_{S'}.$$

(2) Let us suppose that $\Theta_{S,T}(0) \in \mathbf{Z}[G]_{0,S} \cdot \operatorname{Fitt}(A_S)$. One can therefore write $\Theta_{S,T}(0)$ as a finite sum $\Theta_{S,T}(0) = \sum_i \alpha_i \cdot f_i$, where $\alpha_i \in \mathbf{Z}[G]_{0,S}$, and $f_i \in \operatorname{Fitt}(A_S)$. Property (e) above together with §1.4(g) imply that there exist elements $u_i \in \overset{r'}{\wedge} U_{S'}$ so that $f_i = \Phi'(u_i)$, for every index $i$, and therefore (f) implies that $\varepsilon_S = \Theta_{S,T}(0) = \Phi(\sum \alpha_i \cdot u_i)$. The obvious fact that, in our context, $\mathbf{Z}[G]_{0,S} = \mathbf{Z}[G]_{r',S'}$ therefore implies that:

$$\sum_i \alpha_i \cdot u_i \in \mathbf{Z}[G]_{r',S'} \cdot \overset{r'}{\wedge} U_{S'} \subseteq \Lambda_{S'}.$$

Properties (a) and (d) imply then that $\varepsilon_{S'} = \sum \alpha_i \cdot u_i \in \mathbf{Z}[G]_{r',S'} \cdot \overset{r'}{\wedge} U_{S'}$. □

We conclude this section with the following statement, whose proof can be found in [14], §5 (see Corollary 5.4). For a more detailed proof, see [12], Chapter II.

**Proposition 2.3.** *Let us suppose that the set of data $(K/k, S, T, r)$ satisfies hypotheses (H), and that $\varepsilon_{S,T}$ is the element satisfying $\mathbf{Q} \cdot \operatorname{St}(K/k, S, T, r)$. If $\varepsilon_{S,T} \in \mathbf{Z}[1/g]\Lambda_{S,T}$, then the following statements are equivalent:*

(1) $\varepsilon_{S,T} \in \mathbf{Z}[1/g]\operatorname{Fitt}(A_{S,T})\Lambda_{S,T}$.
(2) $\mathbf{Z}[1/g][G]\varepsilon_{S,T} = \mathbf{Z}[1/g]\operatorname{Fitt}(A_{S,T})\Lambda_{S,T}$.

We will actually prove that statement (1) in Proposition 2.3 is always satisfied in the function field case (see Theorem 3.2.1(2)), and therefore (2) is always satisfied (see Corollary 3.2.2). This fact turns out to be of crucial importance in proving that Gras–type Conjectures hold true in the function field setting (see [13]).

## 3. Conjecture B up to primes dividing $|G(K/k)|$

3.1. THE $r = 0$ CASE

In the context described in §1.7, let us suppose that the set of data $(K/k, S, T, 0)$ satisfies hypotheses (H). The goal of this section is the proof of the following:

**Theorem 3.1.1.** *If the set of data $(K/k, S, T, 0)$ satisfies hypotheses* (H)*, then*

$$\Theta_{S,T}(0) \in \mathbf{Z}[1/g] \operatorname{Fitt}_{\mathbf{Z}[1/g][G]}(A_{S,T}).$$

**Remark 1.** Before starting to prove the theorem above, let us first notice that, according to §1.4 (c), its statement is equivalent to

$$\Theta_{S,T}(0) \in \operatorname{Fitt}_{\mathbf{Z}_\ell[G]}(A_{S,T} \otimes \mathbf{Z}_\ell),$$

for all primes $\ell$, such that $\gcd(\ell, g) = 1$. We will therefore prove the Theorem "prime by prime".

Let $S_{X_1}$ and $S_X$ be the finite sets of closed points on $X_1$ and $X$ respectively, lying above points in $Y_1$ associated to primes in $S$. $T_{X_1}$ and $T_X$ have similar meanings. If $Z$ is a subset of closed points of $X_1$ or $X$, then $\operatorname{Div}(Z)$ and $\operatorname{Div}^0(Z)$ will denote the groups of divisors and respectively divisors of degree 0, supported on $Z$. (If $Z \subseteq X$, degree 0 means *multidegree* 0, i.e. degree 0 on each connected component $X^{(i)}$.) Let

$$\operatorname{Pic}^0(X_1) = \frac{\operatorname{Div}^0(X_1)}{\{\operatorname{div}(f) : f \in K^\times\}}, \quad \operatorname{Pic}^0(X) = \frac{\operatorname{Div}^0(X)}{\{\operatorname{div}(f) : f \in \bigoplus_{0 \leq i \leq \nu-1} K^{(i)\times}\}}$$

be the Picard groups associated to $X_1$ and $X$ respectively. As explained in §1.7.3, $\operatorname{Pic}^0(X)$ is the underlying group of an abelian variety of dimension $g_X$, on which the geometric Frobenius map induces a bijective endomorphism $F_{*,1}$. $\operatorname{Pic}^0(X_1)$ can be naturaly viewed as a finite subgroup of $\operatorname{Pic}^0(X)$, and it is precisely the part of $\operatorname{Pic}^0(X)$ fixed by the action of $F_{*,1}$.

We will also consider the following groups:

$$\operatorname{Pic}^0(X_1)_T = \frac{\operatorname{Div}^0(X_1 \setminus T_{X_1})}{\{\operatorname{div}(f) : f \in K^\times, f \equiv 1 \mod w, \forall w \in T_{X_1}\}},$$

$$\operatorname{Pic}^0(X)_T = \frac{\operatorname{Div}^0(X \setminus T_X)}{\{\operatorname{div}(f) : f \in \bigoplus_i K^{(i)\times}, f \equiv 1 \mod w, \forall w \in T_X\}}.$$

It can be shown that $\operatorname{Pic}^0(X)_T$ is the underlying group of an algebraic group as well. We will not use this fact in what follows. $\operatorname{Pic}^0(X_1)_T$ sits inside the part of $\operatorname{Pic}^0(X)_T$ fixed by the action $F_{*,1}$ of the geometric Frobenius morphism and, as Lemma 3.1.4 below shows, it is in fact equal to this. There obviously are surjective group homomorphisms

$$\pi_{X_1} : \operatorname{Pic}^0(X_1)_T \longrightarrow \operatorname{Pic}^0(X_1), \quad \pi_X : \operatorname{Pic}^0(X)_T \longrightarrow \operatorname{Pic}^0(X).$$

There is a commutative diagram of $\mathbf{Z}[G]$–modules, with exact rows

$$
\begin{array}{ccccccccc}
0 & \longrightarrow & K^\times/\mathbf{F}_{q^\nu}^\times & \xrightarrow{\mathrm{div}} & \mathrm{Div}^0(X_1) & \longrightarrow & \mathrm{Pic}^0(X_1) & \longrightarrow & 0 \\
& & \pi_K \downarrow & & \phi_S \downarrow & & \psi_S \downarrow & & \\
0 & \longrightarrow & K^\times/U_S & \xrightarrow{\mathrm{div}_S} & \mathrm{Div}(X_1 \setminus S_{X_1}) & \longrightarrow & A_S & \longrightarrow & 0,
\end{array}
\quad (11)
$$

where:
- $\pi_K$ is the usual projection ($\mathbf{F}_{q^\nu}^\times \subseteq U_S$).
- $\phi_S$ is the "forgetful" map defined by $\phi\left(\sum_{w_i \in X_1} n_i w_i\right) = \sum_{w_i \in X_1 \setminus S_{X_1}} n_i w_i$.
- $\psi_S$ is the map induced by $\phi_S$.
- div is the usual divisor map.
- $\mathrm{div}_S$ is the "forgetful" divisor map, defined by $\mathrm{div}_S(f) = \sum_{w \in X_1 \setminus S_{X_1}} \mathrm{ord}_w(f) \cdot w$, for any $f \in K$.

The snake lemma applied to (11) and the surjectivity of $\pi_K$ give an isomorphism of $\mathbf{Z}[G]$–modules $\mathrm{coker}(\phi_S) \xrightarrow{\sim} \mathrm{coker}(\psi_S)$. Let $d\mathbf{Z} = \mathrm{Im}\left(\mathrm{Div}(S_{X_1}) \xrightarrow{\mathrm{deg}} \mathbf{Z}\right)$, where "deg" is the divisor degree map on $X_1$. It is an easy observation that $\mathrm{coker}(\phi_S) \xrightarrow{\sim} \mathbf{Z}/d\mathbf{Z}$ as $\mathbf{Z}[G]$–modules, with $G$ acting trivially on $\mathbf{Z}/d\mathbf{Z}$. Therefore we obtain the following isomorphism of
$\mathbf{Z}[G]$–modules

$$\mathrm{coker}(\psi_S) \xrightarrow{\sim} \mathbf{Z}/d\mathbf{Z}, \quad (12)$$

with $G$ acting trivially on $\mathbf{Z}/d\mathbf{Z}$.

There is a commutative diagram of $\mathbf{Z}[G]$–modules with exact rows

$$
\begin{array}{ccccccccc}
0 & \longrightarrow & (\bigoplus_{w \in T_{X_1}} \mathbf{F}_{q^\nu}(w)^\times)/\mathbf{F}_{q^\nu}^\times & \xrightarrow{i_T} & \mathrm{Pic}^0(X_1)_T & \xrightarrow{\pi_{X_1}} & \mathrm{Pic}^0(X_1) & \longrightarrow & 0 \\
& & \pi_T \downarrow & & \psi_{S,T} \downarrow & & \psi_S \downarrow & & \\
0 & \longrightarrow & (\bigoplus_{w \in T_{X_1}} \mathbf{F}_{q^\nu}(w)^\times)/U_S & \xrightarrow{i_{S,T}} & A_{S,T} & \longrightarrow & A_S & \longrightarrow & 0,
\end{array}
$$
(13)

where:
- The lower row comes from the exact sequence (1).
- $i_T(\widehat{(a_w)_{w \in T_{X_1}}}) = \widehat{\mathrm{div}(f)}$, for every $(a_w)_{w \in T_{X_1}} \in \bigoplus_{w \in T_{X_1}} \mathbf{F}_{q^\nu}(w)^\times$ and $f \in K$,

satisfying $f \equiv a_w \mod w, \forall w \in T_{X_1}$.
- $\psi_{S,T}(\widehat{\sum_i n_i \cdot w_i}) = $ the ideal class corresponding to $\prod_i w_i^{n_i}$ in $A_{S,T}$.

The snake lemma applied to (13) and the surjectivity of $\pi_T$ imply that $\mathrm{coker}(\psi_{S,T}) \xrightarrow{\sim} \mathrm{coker}(\psi_S)$. (12) therefore gives an isomorphism of $\mathbf{Z}[G]$–modules

$$\mathrm{coker}(\psi_{S,T}) \xrightarrow{\sim} \mathbf{Z}/d\mathbf{Z}, \quad (14)$$

with $G$ acting trivially on $\mathbf{Z}/d\mathbf{Z}$.

**Proposition 3.1.2.** *If $\ell$ is a prime number, such that $\gcd(\ell, g) = 1$, then*

$$\det\nolimits_{\mathbf{Z}_\ell[G]}(1 - F_{*,0}|\mathrm{T}_\ell(\ker \alpha_S)) \in \mathrm{Fitt}_{\mathbf{Z}_\ell[G]}(\mathrm{coker}(\psi_{S,T}) \otimes \mathbf{Z}_\ell).$$

Before starting to prove the proposition, let us observe that the determinant in the statement makes sense due to the $\mathbf{Z}_\ell[G]$–projectivity of $\mathrm{T}_\ell(\ker \alpha_S)$ (see Remark 2, §1.7.4 and Remark 1, §1.3). In what follows, if $M$ is an abelian group and $f \in \mathrm{End}(M)$, then $M^f$ denotes the subgroup of $M$ fixed by $f$.

**Proof.** Since $\alpha_S$ is $F_{*,0}$–equivariant, we have a commutative diagram of $\mathbf{Z}[G]$–modules, with exact rows

$$\begin{array}{ccccccccc}
0 & \longrightarrow & \ker \alpha_S & \longrightarrow & \mathrm{Div}(S_X) & \xrightarrow{\alpha_S} & \mathrm{H}_0(X, \mathbf{Z}) & \longrightarrow & 0 \\
& & {\scriptstyle (1-F_{*,0})_1} \downarrow & & {\scriptstyle (1-F_{*,0})_2} \downarrow & & {\scriptstyle (1-F_{*,0})_3} \downarrow & & \\
0 & \longrightarrow & \ker \alpha_S & \longrightarrow & \mathrm{Div}(S_X) & \longrightarrow & \mathrm{H}_0(X, \mathbf{Z}) & \longrightarrow & 0.
\end{array}$$

Obviously, we have the following equalities:
- $\ker(1 - F_{*,0})_2 = \mathrm{Div}(S_X)^{F_{*,0}} = \mathrm{Div}(S_{X_1})$,
- $\ker(1 - F_{*,0})_3 =$ the rank 1, free $\mathbf{Z}$–module generated by $A = \sum\limits_{0 \leq i \leq \nu - 1} [X^{(i)}]$,

with trivial $G$–action.

The snake lemma gives therefore a long exact sequence of $\mathbf{Z}[G]$–modules

$$0 \longrightarrow \ker(1 - F_{*,0})_1 \longrightarrow \mathrm{Div}(S_{X_1}) \xrightarrow{A \cdot \deg} \mathbf{Z} \cdot A \longrightarrow$$
$$\longrightarrow \mathrm{coker}(1 - F_{*,0})_1 \longrightarrow \mathrm{coker}(1 - F_{*,0})_2 \longrightarrow \mathrm{coker}(1 - F_{*,0})_3 \longrightarrow 0$$

Since $\mathrm{Im}(A \cdot \deg) = d\mathbf{Z} \cdot A$, and $G$ fixes $A$, we get an injective morphism of $\mathbf{Z}[G]$–modules

$$\mathbf{Z}/d\mathbf{Z} \hookrightarrow \mathrm{coker}(1 - F_{*,0})_1, \tag{15}$$

with $G$ acting trivially on $\mathbf{Z}/d\mathbf{Z}$. Since for $\ell \nmid g$, $\mathbf{Z}_\ell[G]$ is a direct sum of Dedekind domains, (14), (15) and (f), §1.4 imply that

$$\mathrm{Fitt}_{\mathbf{Z}_\ell[G]}(\mathrm{coker}(1 - F_{*,0})_1 \otimes \mathbf{Z}_\ell) \subseteq \mathrm{Fitt}_{\mathbf{Z}_\ell[G]}(\mathrm{coker}(\psi_{S,T}) \otimes \mathbf{Z}_\ell), \quad \forall \ell, \ell \nmid g. \tag{16}$$

If one considers the decomposition $\mathbf{Z}_\ell[G] = \bigoplus\limits_{\psi \in \widehat{G}(\mathbf{Q}_\ell)} D_\psi$, as in §1.3, then the fact that $\mathrm{T}_\ell(\ker \alpha_S)$ is a projective $\mathbf{Z}_\ell[G]$–module implies that one has an isomorphism $\mathrm{T}_\ell(\ker \alpha_S) \xrightarrow{\sim} \oplus_\psi D_\psi^{n_\psi}$, for some positive integers $n_\psi$ (see Remark 1, §1.3). The exact sequence

$$\mathrm{T}_\ell(\ker \alpha_S) \xrightarrow{(1-F_{*,0})_1} \mathrm{T}_\ell(\ker \alpha_S) \longrightarrow \mathrm{coker}(1 - F_{*,0})_1 \otimes \mathbf{Z}_\ell \longrightarrow 0$$

gives therefore a situation similar to the one described in §1.4 (d), with $R = \mathbf{Z}_\ell[G]$ and $N = \mathrm{T}_\ell(\ker \alpha_S)$. We therefore have

$$\det\nolimits_{\mathbf{Z}_\ell[G]}(1 - F_{*,0}|\mathrm{T}_\ell(\ker \alpha_S)) \in \mathrm{Fitt}_{\mathbf{Z}_\ell[G]}(\mathrm{coker}(1 - F_{*,0})_1 \otimes \mathbf{Z}_\ell). \tag{17}$$

Relations (16) and (17) imply the statement in Proposition 3.1.2. □

**Proposition 3.1.3.** *If $\ell$ is a prime number, such that $\gcd(\ell, g) = 1$, then*

$$\det_{\mathbf{Z}_\ell[G]} (1 - F_{*,2}|\, \mathrm{T}_\ell(\mathrm{coker}(\beta_T))) \cdot \det_{\mathbf{Z}_\ell[G]} (1 - F_{*,1}|\, \mathrm{H}_1(X, \mathbf{Z}_\ell)) \in$$
$$\in \mathrm{Fitt}_{\mathbf{Z}_\ell[G]} \left( \mathrm{Pic}^0(X_1)_T \otimes \mathbf{Z}_\ell \right).$$

**Proof.** Let us recall that we have an exact sequence of $\mathbf{Z}[G]$–modules

$$0 \longrightarrow \mathrm{H}_0(X, \mathbf{Z}) \otimes \mathbf{F}^\times \xrightarrow{\beta_T} \mathrm{H}_0(T_X, \mathbf{Z}) \otimes \mathbf{F}^\times \longrightarrow \mathrm{coker}(\beta_T) \longrightarrow 0. \qquad (18)$$

The $q$–power geometric Frobenius morphism induces isomorphisms $F_{*,2} = F_{*,0} \otimes \sigma_q$ on each term of (18), where $\sigma_q : \mathbf{F}^\times \longrightarrow \mathbf{F}^\times$ is the $q$–power map, and $\beta_T$ is $F_{*,2}$–invariant (see §1.7.4).

Let $\Sigma = G(\mathbf{F}/\mathbf{F}_q)$ (topologically generated by $\sigma_q$), let $\Sigma_\nu$ be the closed subgroup of $\Sigma$ (topologically) generated by $\sigma_q^\nu$, and let $\Gamma_0$ be the finite subgroup of $\mathrm{Aut}_{\mathbf{Z}}(\mathrm{H}_0(T_X, \mathbf{Z}))$ generated by $F_{*,0}$. Then obviously the profinite group $\Gamma_0 \times \Sigma$ acts continously on each term of (18), considered with the discrete topology.

Let $H$ be the closed subgroup of $\Gamma_0 \times \Sigma$, (topologically) generated by $F_{*,2} = F_{*,0} \otimes \sigma_q$. Let $H_\nu$ be the closed subgroup of $H$ generated by $F_{*,2}^\nu = F_{*,0}^\nu \otimes \sigma_q^\nu$. Since the actions of $H$ and $G$ on $\mathrm{H}_0(T_X, \mathbf{Z}) \otimes \mathbf{F}^\times$ commute, (18) is an exact sequence of $H \times G$–modules.

We also have an exact sequence of $G$–modules

$$0 \longrightarrow \mathrm{coker}(\beta_T) \longrightarrow \mathrm{Pic}^0(X)_T \longrightarrow \mathrm{Pic}^0(X) \longrightarrow 0, \qquad (19)$$

which preserves the Frobenius action on each of its terms (see [17], Chpt.V).

Let $\ell$ be a prime number. Remark 3, §1.7.4 shows that the multiplication by $\ell$ maps

$$\mathrm{coker}(\beta_T)\left[\ell^{n+1}\right] \xrightarrow{\ell} \mathrm{coker}(\beta_T)\left[\ell^n\right]$$

are surjective for every $n \geq 0$. The definition of $\mathrm{Pic}^0(X_1)_T$, combined with the weak approximation theorem for valuations on $K^{(i)}$, for all $0 \leq i \leq \nu - 1$, easily implies that the natural maps

$$\mathrm{Pic}^0(X_1)_T[\ell^n] \longrightarrow \mathrm{Pic}^0(X_1)[\ell^n]$$

are surjective as well, for every $n \geq 0$. These two facts show that the short exact sequence (19) gives a short exact sequence of $G$–modules

$$0 \longrightarrow \mathrm{T}_\ell(\mathrm{coker}(\beta_T)) \longrightarrow \mathrm{T}_\ell\left(\mathrm{Pic}^0(X)_T\right) \longrightarrow \mathrm{T}_\ell\left(\mathrm{Pic}^0(X)\right) \longrightarrow 0, \qquad (20)$$

at the level of $\ell$–adic Tate modules, preserving the Frobenius action on its terms.

We will need two lemmas:

**Lemma 3.1.4.** *We have the following:*
   (1) $\mathrm{T}_\ell\left(\mathrm{Pic}^0(X)_T\right)$ *is a free* $\mathbf{Z}_\ell$*-module of finite rank, for every prime* $\ell$.
   (2) $\mathrm{Pic}^0(X)_T^{F_{*,1}} \xrightarrow{\sim} \mathrm{Pic}^0(X_1)_T$ *as* $\mathbf{Z}[G]$*-modules, where* $\mathrm{Pic}^0(X)_T^{F_{*,1}}$ *is the part of* $\mathrm{Pic}^0(X)_T$ *fixed by* $F_{*,1}$.

**Proof.** Statement (1) follows directly from the exact sequence (20) and from the fact that both $\mathrm{T}_\ell(\mathrm{coker}(\beta_T))$ and $\mathrm{T}_\ell\left(\mathrm{Pic}^0(X)\right)$ are free $\mathbf{Z}_\ell$–modules of finite rank (see §1.7.3 and Remark3, §1.7.4).

(2) The fact that $F_{*,0}^\nu = 1_{\mathrm{H}_0(X,\mathbf{Z})}$ implies that there is an $H$–isomorphism

$$\mathrm{H}_0(X,\mathbf{Z}) \otimes \mathbf{F}^\times \xrightarrow{\sim} \left(\mathbf{Z} \cdot \left[X^{(0)}\right] \otimes \mathbf{F}^\times\right) \bigotimes_{\mathbf{Z}[H_\nu]} \mathbf{Z}[H]. \tag{21}$$

Shapiro's Lemma and Hilbert's Theorem 90 imply therefore that, at the level of Galois cohomology groups, we have

$$\mathrm{H}^i\left(H, \mathrm{H}_0(X,\mathbf{Z}) \otimes \mathbf{F}^\times\right) \xrightarrow{\sim} \mathrm{H}^i\left(H_\nu, \mathbf{Z} \cdot \left[X^{(0)}\right] \otimes \mathbf{F}^\times\right) = \mathrm{H}^i\left(\Sigma_\nu, \mathbf{F}^\times\right) = 0, \tag{22}$$

for $i = 1, 2$.

The fact that $H/H_\nu$ is finite cyclic and (21) above imply that there are isomorphisms of $G$–modules

$$\left(\mathrm{H}_0(X,\mathbf{Z}) \otimes \mathbf{F}^\times\right)^H \xrightarrow{\sim} \left(\mathbf{Z} \cdot \left[X^{(0)}\right] \otimes \mathbf{F}^\times\right)^{H_\nu} \xrightarrow{\sim} \mathbf{F}^{\times \Sigma_\nu} = \mathbf{F}_{q^\nu}^\times. \tag{23}$$

In a similar fashion one can prove that

$$\mathrm{H}^1\left(H, \mathrm{H}_0(T_X, \mathbf{Z}) \otimes \mathbf{F}^\times\right) = 0 \tag{24}$$

and that there is an isomorphism of $G$–modules

$$\left(\mathrm{H}_0(T_X, \mathbf{Z}) \otimes \mathbf{F}^\times\right)^H \xrightarrow{\sim} \bigoplus_{w \in T_X} \mathbf{F}_{q^\nu}(w)^\times. \tag{25}$$

As a consequence of (22–25), the long exact sequence of $H$–cohomology groups associated to the short exact sequence (18) starts as

$$0 \longrightarrow \mathbf{F}_{q^\nu}^\times \longrightarrow \bigoplus_{w \in T_{X_1}} \mathbf{F}_{q^\nu}(w)^\times \longrightarrow (\mathrm{coker}(\beta_T))^H \longrightarrow$$

$$\longrightarrow 0 \longrightarrow 0 \longrightarrow \mathrm{H}^1(H, \mathrm{coker}(\beta_T)) \longrightarrow 0 \longrightarrow \cdots.$$

This implies that we have the following $G$–isomorphisms:

$$\mathrm{H}^1(H, \mathrm{coker}(\beta_T)) = 0 \quad \text{and} \quad (\mathrm{coker}(\beta_T))^H \xrightarrow{\sim} (\bigoplus_{w \in T_{X_1}} \mathbf{F}_{q^\nu}(w)^\times)/\mathbf{F}_{q^\nu}^\times. \tag{26}$$

Let $\mathcal{H}$ be the subgroup of $\mathrm{Aut}_{\mathbf{Z}}\left(\mathrm{Pic}^0\left(X_1\right)_T\right)$ generated by $F_{*,1}$. As we mentioned above, the restriction of $F_{*,1}$ to $\mathrm{coker}\,(\beta_T)$ via the injective morphism $\mathrm{coker}\,(\beta_T) \longrightarrow \mathrm{Pic}^0\left(X_1\right)_T$ coincides to $F_{*,2}$. This observation together with (26) show that the short exact sequence (19) gives a long exact sequence of $\mathcal{H}$–cohomology groups which starts as

$$0 \longrightarrow (\mathrm{coker}(\beta_T))^{F_{*,2}} \longrightarrow \mathrm{Pic}^0(X)_T^{F_{*,1}} \longrightarrow \mathrm{Pic}^0(X)^{F_{*,1}} \longrightarrow 0.$$

We therefore have a commutative diagram with exact rows

$$\begin{array}{ccccccccc}
0 & \longrightarrow & (\bigoplus_{w \in T_{X_1}} \mathbf{F}_{q^\nu}(w)^\times)/\mathbf{F}_{q^\nu}^\times & \xrightarrow{i_T} & \mathrm{Pic}^0\left(X_1\right)_T & \xrightarrow{\pi_{X_1}} & \mathrm{Pic}^0\left(X_1\right) & \longrightarrow & 0 \\
& & \downarrow{i_1} & & \downarrow{i_2} & & \downarrow{i_3} & & \\
0 & \longrightarrow & (\mathrm{coker}(\beta_T))^{F_{*,2}} & \longrightarrow & \mathrm{Pic}^0(X)_T^{F_{*,1}} & \longrightarrow & \mathrm{Pic}^0(X)^{F_{*,1}} & \longrightarrow & 0,
\end{array}$$
(27)

where $i_2$ and $i_3$ are the natural inclusions and $i_1$ is the isomorphism in (26). We already know that $i_3$ and $i_1$ are isomorphisms (see comments at the beginning of this section and (26) above, respectively). The snake lemma therefore implies that $i_2$ is an isomorphism as well, and this concludes the proof of (2). □

**Lemma 3.1.5.** *Let $M$ be a $\mathbf{Z}[G]$–module, which has nontrivial $\mathbf{Z}$–torsion, let $f \in End_{\mathbf{Z}[G]}(M)$, and $\ell$ a prime number, such that*

(1) $\mathrm{T}_\ell(M)$ *is a free $\mathbf{Z}_\ell$–module of finite rank.*
(2) $M^f \stackrel{\mathrm{def}}{=} \{m \in M : f(m) = m\}$ *is finite.*

*Then there is an isomorphism of $\mathbf{Z}_\ell[G]$–modules*

$$\mathrm{T}_\ell(M)/(1-f)\mathrm{T}_\ell(M) \xrightarrow{\sim} \mathbf{Z}_\ell \otimes M^f.$$

**Proof.** Let us consider the exact sequence of $\mathbf{Z}_\ell$–modules

$$0 \longrightarrow \mathbf{Z}_\ell \longrightarrow \mathbf{Q}_\ell \longrightarrow \mathbf{Q}_\ell/\mathbf{Z}_\ell \longrightarrow 0.$$

In virtue of hypothesis (1) above,

$$0 \longrightarrow \mathrm{T}_\ell(M) \longrightarrow \mathrm{T}_\ell(M) \bigotimes_{\mathbf{Z}_\ell} \mathbf{Q}_\ell \longrightarrow \mathrm{T}_\ell(M) \bigotimes_{\mathbf{Z}_\ell} \mathbf{Q}_\ell/\mathbf{Z}_\ell \longrightarrow 0$$

is still exact. We therefore get a commutative diagram with exact rows

$$\begin{array}{ccccccccc}
0 & \longrightarrow & \mathrm{T}_\ell(M) & \longrightarrow & \mathrm{T}_\ell(M) \bigotimes_{\mathbf{Z}_\ell} \mathbf{Q}_\ell & \longrightarrow & \mathrm{T}_\ell(M) \bigotimes_{\mathbf{Z}_\ell} \mathbf{Q}_\ell/\mathbf{Z}_\ell & \longrightarrow & 0 \\
& & \downarrow{(1-f)_1} & & \downarrow{(1-f)_2} & & \downarrow{(1-f)_3} & & \\
0 & \longrightarrow & \mathrm{T}_\ell(M) & \longrightarrow & \mathrm{T}_\ell(M) \bigotimes_{\mathbf{Z}_\ell} \mathbf{Q}_\ell & \longrightarrow & \mathrm{T}_\ell(M) \bigotimes_{\mathbf{Z}_\ell} \mathbf{Q}_\ell/\mathbf{Z}_\ell & \longrightarrow & 0,
\end{array}$$

where $(1-f)_1 = (1-f)$ at the level of Tate modules, $(1-f)_2 = (1-f)_1 \otimes 1_{\mathbf{Q}_\ell}$, and $(1-f)_3 = (1-f)_1 \otimes 1_{\mathbf{Q}_\ell/\mathbf{Z}_\ell}$. Hypothesis (2) above implies that $(1-f)_1$ is injective and therefore (1) implies that $(1-f)_2$ is an isomorphism of $\mathbf{Q}_\ell$–vector spaces. Snake lemma applied to the diagram above gives therefore an isomorphism of $\mathbf{Z}_\ell[G]$–modules

$$\ker(1-f)_3 \xrightarrow{\sim} \mathrm{T}_\ell(M)/(1-f)\mathrm{T}_\ell(M). \tag{28}$$

On the other hand

$$\mathrm{T}_\ell(M) \bigotimes_{\mathbf{Z}_\ell} \mathbf{Q}_\ell/\mathbf{Z}_\ell \xrightarrow{\sim} \bigcup_{i \geq 0} M[\ell^i].$$

Let $M[\ell^\infty] = \bigcup_{i \geq 0} M[\ell^i]$. Then the commutative diagram

$$\begin{array}{ccc} \mathrm{T}_\ell(M) \bigotimes_{\mathbf{Z}_\ell} \mathbf{Q}_\ell/\mathbf{Z}_\ell & \xrightarrow{\sim} & M[\ell^\infty] \\ {\scriptstyle (1-f)_3} \downarrow & & \downarrow {\scriptstyle (1-f)} \\ \mathrm{T}_\ell(M) \bigotimes_{\mathbf{Z}_\ell} \mathbf{Q}_\ell/\mathbf{Z}_\ell & \xrightarrow{\sim} & M[\ell^\infty], \end{array}$$

in addition to (28) and hypothesis (2) give the isomorphisms of $\mathbf{Z}_\ell[G]$–modules

$$\mathrm{T}_\ell(M)/(1-f)\mathrm{T}_\ell(M) \xrightarrow{\sim} M[\ell^\infty]^f \xrightarrow{\sim} M^f \otimes \mathbf{Z}_\ell.$$

This concludes the proof of Lemma 3.1.5. □

We now return to the proof of Proposition 3.1.3. Lemma 3.1.4 and the finiteness of $\mathrm{Pic}^0(X_1)_T$ (see the upper row of commutative diagram (27) and keep in mind that $\mathrm{Pic}^0(X_1)$ is finite) imply that $\mathrm{Pic}^0(X_1)_T$ together with its endomorphism $F_{*,1}$, and any prime number $\ell$, satisfy the hypotheses in Lemma 3.1.5. We therefore have an isomorphism of $\mathbf{Z}_\ell[G]$–modules

$$\mathrm{T}_\ell\left(\mathrm{Pic}^0(X)_T\right)/(1-F_{*,1})\mathrm{T}_\ell\left(\mathrm{Pic}^0(X)_T\right) \xrightarrow{\sim} \mathrm{Pic}^0(X_1)_T \otimes \mathbf{Z}_\ell, \tag{29}$$

for any prime number $\ell$. The exact sequence (20) gives a commutative diagram of $\mathbf{Z}_\ell[G]$–modules with exact rows

$$\begin{array}{ccccccccc} 0 & \longrightarrow & \mathrm{T}_\ell(\mathrm{coker}(\beta_T)) & \longrightarrow & \mathrm{T}_\ell\left(\mathrm{Pic}^0(X)_T\right) & \longrightarrow & \mathrm{T}_\ell\left(\mathrm{Pic}^0(X)\right) & \longrightarrow & 0 \\ & & {\scriptstyle (1-F_{*,2})} \downarrow & & {\scriptstyle (1-F_{*,1})_1} \downarrow & & {\scriptstyle (1-F_{*,1})_2} \downarrow & & \\ 0 & \longrightarrow & \mathrm{T}_\ell(\mathrm{coker}(\beta_T)) & \longrightarrow & \mathrm{T}_\ell\left(\mathrm{Pic}^0(X)_T\right) & \longrightarrow & \mathrm{T}_\ell\left(\mathrm{Pic}^0(X)\right) & \longrightarrow & 0 \end{array}$$

The finiteness of $\mathrm{Pic}^0(X)_T^{F_{*,1}} = \mathrm{Pic}^0(X_1)_T$ and $\mathrm{Pic}^0(X)^{F_{*,1}} = \mathrm{Pic}^0(X_1)$ respectively implies that the maps $(1-F_{*,1})_1$ and $(1-F_{*,1})_2$ are both injective. Snake lemma applied to the diagram above, and (29), give an exact sequence of $\mathbf{Z}_\ell[G]$–modules

$$0 \longrightarrow \mathrm{coker}(1 - F_{*,2}) \longrightarrow \mathrm{Pic}^0(X_1)_T \otimes \mathbf{Z}_\ell \longrightarrow \mathrm{coker}(1 - F_{*,1})_2 \longrightarrow 0$$

Property (e) §1.4 implies that

$$\begin{aligned}\mathrm{Fitt}_{\mathbf{Z}_\ell[G]}(\mathrm{coker}(1 - F_{*,1})_2) \cdot \mathrm{Fitt}_{\mathbf{Z}_\ell[G]}(\mathrm{coker}(1 - F_{*,2})) &\subseteq \\ \in \mathrm{Fitt}_{\mathbf{Z}_\ell[G]}\left(\mathrm{Pic}^0(X_1)_T \otimes \mathbf{Z}_\ell\right).&\end{aligned} \quad (30)$$

On the other hand, both $\mathrm{T}_\ell(\mathrm{coker}(\beta_T))$ and $\mathrm{T}_\ell\left(\mathrm{Pic}^0(X)\right)$ are free $\mathbf{Z}_\ell$–modules, and therefore projective $\mathbf{Z}_\ell[G]$–modules (see Remark 1 §1.3). An argument similar to the one used in proving (17) shows that

$$\begin{aligned}\mathrm{det}_{\mathbf{Z}_\ell[G]}(1 - F_{*,2}|\mathrm{T}_\ell(\mathrm{coker}(\beta_T))) &\in \mathrm{Fitt}_{\mathbf{Z}_\ell[G]}(\mathrm{coker}(1 - F_{*,2})) \quad \text{and} \\ \mathrm{det}_{\mathbf{Z}_\ell[G]}((1 - F_{*,1})_1|\mathrm{T}_\ell(\mathrm{Pic}^0(X))) &\in \mathrm{Fitt}_{\mathbf{Z}_\ell[G]}(\mathrm{coker}(1 - F_{*,1})_2).\end{aligned} \quad (31)$$

According to (30) and (31) we can finally conclude that Proposition 3.1.3 is true. $\square$

**Proof of Theorem 3.1.1.** Let us consider a prime $\ell$, such that $\gcd(\ell, g) = 1$, and the exact sequence of $\mathbf{Z}_\ell[G]$–modules

$$\mathrm{Pic}^0(X_1)_T \otimes \mathbf{Z}_\ell \xrightarrow{\psi_{S,T} \otimes 1_{\mathbf{Z}_\ell}} A_{S,T} \otimes \mathbf{Z}_\ell \longrightarrow \mathrm{coker}(\psi_{S,T}) \otimes \mathbf{Z}_\ell \longrightarrow 0$$

(see diagram (21) for the definition of $\psi_{S,T}$). Property (e) §1.4 implies that

$$\begin{aligned}\mathrm{Fitt}_{\mathbf{Z}_\ell[G]}(\mathrm{coker}(\psi_{S,T}) \otimes \mathbf{Z}_\ell) \cdot \mathrm{Fitt}_{\mathbf{Z}_\ell[G]}(\mathrm{Pic}^0(X_1)_T \otimes \mathbf{Z}_\ell) &\subseteq \\ \subseteq \mathrm{Fitt}_{\mathbf{Z}_\ell[G]}(A_{S,T} \otimes \mathbf{Z}_\ell).&\end{aligned} \quad (32)$$

**Step 1.** Let us first suppose that $\ell \neq p$. Then Propositions 3.1.2 and 3.1.3, in addition to (32) above, and Lemma 1.7.4.4 (1) for $s = 0$ (and therefore $u = 1$), imply that

$$\begin{aligned}\Theta_{S,T}(0) = &\mathrm{det}_{\mathbf{Z}_\ell[G]}(1 - F_{*,0}|\mathrm{T}_\ell(\ker \alpha_S)) \cdot \mathrm{det}_{\mathbf{Z}_\ell[G]}(1 - F_{*,1}|\mathrm{H}_1(X, \mathbf{Z}_\ell)) \cdot \\ &\cdot \mathrm{det}_{\mathbf{Z}_\ell[G]}(1 - F_{*,2}|\mathrm{T}_\ell(\mathrm{coker}(\beta_T))) \in \mathrm{Fitt}_{\mathbf{Z}_\ell[G]}(A_{S,T} \otimes \mathbf{Z}_\ell).\end{aligned}$$

**Step 2.** Let us suppose now that $\ell = p$, and $\gcd(p, g) = 1$. The first observation we have to make is that $\mathrm{T}_p(\mathrm{coker}(\beta_T)) = 1$ and therefore Proposition 3.1.3 implies that

$$\mathrm{det}_{\mathbf{Z}_p[G]}(1 - F_{*,1}|\mathrm{H}_1(X, \mathbf{Z}_p)) \in \mathrm{Fitt}_{\mathbf{Z}_p[G]}\left(\mathrm{Pic}^0(X_1)_T \otimes \mathbf{Z}_p\right). \quad (33)$$

Proposition 3.1.2, (32) and (33) above, in addition to Lemma 1.7.4.4 (2) imply therefore that

$$\begin{aligned}\Theta_{S,T}(0) = &Q(1) \cdot \mathrm{det}_{\mathbf{Z}_p[G]}(1 - F_{*,0}|\mathrm{T}_p(\ker \alpha_S)) \cdot \\ &\cdot \mathrm{det}_{\mathbf{Z}_p[G]}((1 - F_{*,1})_1|\mathrm{H}_1(X, \mathbf{Z}_p)) \in \mathrm{Fitt}_{\mathbf{Z}_p[G]}(A_{S,T} \otimes \mathbf{Z}_p).\end{aligned}$$

According to Remark 1, Step 1 and Step 2 above conclude the proof of Theorem 3.1.1. $\square$

## 3.2. ARBITRARY $r$

The following theorem shows that, for general abelian extensios of function fields, Conjecture B is true, up to primes dividing the order of the Galois group.

**Theorem 3.2.1.** *If the set of data $(K/k, S, T, r)$ satisfies hypotheses* (H), *then*

(1) *There exists a unique element $\varepsilon_{S,T} \in \mathbf{Z}[1/g]\Lambda_{S,T}$ satisfying*

$$\Theta_{S,T}^{(r)}(0) = R_W(\varepsilon_{S,T}).$$

(2) *The element $\varepsilon_{S,T}$ above satisfies*

$$\varepsilon_{S,T} \in \mathbf{Z}[1/g]\operatorname{Fitt}(A_{S,T})\Lambda_{S,T}.$$

**Proof.** Let $\{v_1, \ldots, v_r\}$ be the distinguished set of $r$ distinct primes in $S$, which split completely in $K/k$, and let $S_0 = S \setminus \{v_1, \cdots, v_r\}$. Then obviously the set of data $(K/k, S_0, T, 0)$ satisfies hypotheses (H). Theorem 3.1.1 implies that

$$\Theta_{S_0,T}(0) \in \mathbf{Z}[1/g]\operatorname{Fitt}(A_{S_0,T}).$$

Corollary 2.2 (1), with $S' = S$ and $S = S_0$, implies therefore that the element $\varepsilon_{S,T}$ satisfying $\mathbf{Q} \cdot \operatorname{St}(K/k, S, T, r)$, also satisfies $\varepsilon_{S,T} \in \mathbf{Z}[1/g]\operatorname{Fitt}(A_{S,T})\Lambda_{S,T}$. This settles the proof of both (1) and (2) above, as the uniqueness of $\varepsilon_{S,T}$ is a consequence of Remark 5, §1.6. □

**Corollary 3.2.2 ( Raw form of Gras' Conjecture ).** *The element $\varepsilon_{S,T}$, whose existence is proved in the Theorem 3.2.1, satisfies*

$$\mathbf{Z}[1/g][G] \cdot \varepsilon_{S,T} = \mathbf{Z}[1/g]\operatorname{Fitt}(A_{S,T})\Lambda_{S,T} = \operatorname{Fitt}(A_{S,T})\left(\mathbf{Z}[1/g]\overset{r}{\wedge}U_{S,T}\right)_{r,S}.$$

**Proof.** This is a direct consequence of Theorem 3.2.1 (2), Proposition 2.3 and Lemma 1.5.3(1). □

In [13] it will become apparent that Corollary 3.2.2 implies Gras-type conjectures for function fields.

## 4. Proof of a stronger form of Conjecture B for constant field extensions

Througout this section $K/k$ is a constant field extension of function fields of characteristic $p$.

## 4.1. PRELIMINARY CONSIDERATIONS

Let $\mathbf{F}_q$ and $\mathbf{F}_{q^\nu}$ be the exact fields of constants of $k$ and $K$ respectively. Then clearly $G(K/k) \xrightarrow{\sim} G(\mathbf{F}_{q^\nu}/\mathbf{F}_q)$ is cyclic of order $\nu$, with a distinguished generator $\sigma$, satisfying $\sigma(\zeta) = \zeta^q$, for every $\zeta$ in $\mathbf{F}_{q^\nu}$, and $K/k$ is everywhere unramified. We remind the reader that for primes $v$ in $k$, and $w$ in $K$, such that $w|v$, $\sigma_v$ denotes the Frobenius morphism of $v$ relative to $K/k$, $d_v = [\mathbf{F}_q(v) : \mathbf{F}_q]$ and $d_w = [\mathbf{F}_{q^\nu}(w) : \mathbf{F}_{q^\nu}]$. Simple arithmetic considerations show that $\sigma_v = \sigma^{d_v}$, that there are precisely $r_v = \gcd(d_v, \nu)$ primes in $K$ lying above $v$, and that $d_w = d_v/r_v$. (This implies in particular that $v$ splits completely in $K/k$ if and only if $\nu | d_v$.)

**Lemma 4.1.1.** *For any prime number $\ell$, there are $\mathbf{Z}_\ell[G]$–isomorphisms*

(1) $\mathrm{H}_0(X, \mathbf{Z}_\ell) \xrightarrow{\sim} \mathbf{Z}_\ell[G]$, *which sends the action of $F_{*,0}$ on $\mathrm{H}_0(X, \mathbf{Z}_\ell)$ into multiplication by $\sigma^{-1}$ on $\mathbf{Z}_\ell[G]$.*

(2) $\mathrm{H}_2(X, \mathbf{Z}_\ell) \xrightarrow{\sim} \begin{cases} \mathbf{Z}_\ell[G], & \text{if } \ell \neq p \\ 0, & \text{if } \ell = p, \end{cases}$
*which sends the action of $F_{*,2}$ on $\mathrm{H}_2(X, \mathbf{Z}_\ell)$ into multiplication by $q\sigma^{-1}$ on $\mathbf{Z}_\ell[G]$.*

(3) $\mathrm{H}_1(X, \mathbf{Z}_\ell) \xrightarrow{\sim} \begin{cases} \mathbf{Z}_\ell[G]^{2g_X}, & \text{if } \ell \neq p, \\ \mathbf{Z}_p[G]^n, & \text{if } \ell = p, \end{cases}$
*for some integer $n$, $0 \leq n < 2g_X$.*

**Proof.** Let us first notice that

$$F_{*,0}\left(\left[X^{(i)}\right]\right) = \sigma_{*,0}^{-1}\left(\left[X^{(i)}\right]\right), \text{ for every } 0 \leq i \leq \nu - 1. \tag{34}$$

(1) We know that $\mathrm{H}_0(X, \mathbf{Z})$ is $\mathbf{Z}$–free of basis $\{[X^{(i)}] : 0 \leq i \leq \nu - 1\}$, that $G$ is cyclic of order $\nu$ and it permutes the elements $[X^{(i)}]$ transitively. We have also chosen the indices $i$ so that $[X^{(i)}] = F_{*,0}^i([X^{(0)}])$. (34) implies therefore that we can define an isomorphism of $\mathbf{Z}[G]$–modules

$$\mathrm{H}_0(X, \mathbf{Z}) \xrightarrow{\rho}{\sim} \mathbf{Z}[G], \tag{35}$$

given by $\rho([X^{(i)}]) = \sigma^{-i}, \forall 0 \leq i \leq \nu - 1$. Obviously, under this isomorphism, the action of $F_{*,0}$ on the left corresponds to multiplication by $\sigma^{-1}$ on the right. Statement (1) in the lemma can now be obtained by tensoring (35) with $\mathbf{Z}_\ell$.

(2) We already know (see §1.7.3) that there is an isomorphism

$$\mathrm{H}_2(X, \mathbf{Z}_\ell) \xrightarrow{f}{\sim} \mathrm{H}_0(X, \mathbf{Z}_\ell)$$

which takes the action of $F_{*,2}$ on the right hand–side into the action of $qF_{*,0}$ on the left hand–side. The map $(\rho \otimes 1_{\mathbf{Z}_\ell}) \circ f$ gives therefore a $\mathbf{Z}_\ell[G]$–isomorphism satisfying the requirements in (2).

(3) The fact that $|G| = \nu$ and that $G$ permutes transitively and isomorphically the $\nu$ irreducible components $X^{(i)}$ of $X$, and therefore their corresponding Picard groups $\text{Pic}^0(X^{(i)})$, implies that there is an isomorphism of $\mathbf{Z}[G]$–modules

$$\text{H}_1(X, \mathbf{Z}) = \bigoplus_{0 \leq i \leq \nu - 1} \text{Pic}^0(X^{(i)}) \xrightarrow{\sim} \mathbf{Z}[G] \bigotimes_{\mathbf{Z}} \text{Pic}^0(X^{(0)}),$$

with $G$ acting trivially on $\text{Pic}^0(X^{(0)})$. Passing to $\ell$-adic Tate modules, one obtains the desired $\mathbf{Z}_\ell[G]$–isomorphisms. $\square$

4.2. THE $r = 0$ CASE

Throughout this section we will assume that the set of data $(K/k, S, T, 0)$ satisfies hypotheses (H).

**Lemma 4.2.1.** *If $u = q^{-s}$, we have the following:*

(1) *For any prime number $\ell \neq p$*

$$\Theta_{S,T}(u) = \frac{\prod\limits_{v \in S}(1 - (\sigma^{-1}u)^{d_v})}{(1 - \sigma^{-1}u)} \cdot \frac{\prod\limits_{v \in T}(1 - (q\sigma^{-1}u)^{d_v})}{(1 - (q\sigma^{-1}u))} \cdot \det\nolimits_{\mathbf{Z}_\ell[G]}(1 - F_{*,1} \cdot u | \text{H}_1(X, \mathbf{Z}_\ell)).$$

(2) *There exists a polynomial $Q(u) \in \mathbf{Z}_p[G][u]$, such that*

$$\Theta_{S,T}(u) = Q(u) \cdot \frac{\prod\limits_{v \in S}(1-(\sigma^{-1}u)^{d_v})}{(1-\sigma^{-1}u)} \cdot \frac{\prod\limits_{v \in T}(1-(q\sigma^{-1}u)^{d_v})}{(1-(q\sigma^{-1}u))} \cdot$$
$$\cdot \det\nolimits_{\mathbf{Z}_p[G]}(1 - F_{*,1} \cdot u | \text{H}_1(X, \mathbf{Z}_p)).$$

**Proof.** (1) Let $\ell \neq p$ be a prime number. The freeness of the $\text{H}_i(X, \mathbf{Z}_\ell)$'s as $\mathbf{Z}_\ell[G]$–modules implied by Lemma 4.1.1, and Theorem 1.7.4.2, show that

$$\Theta(s) = \prod_{0 \leq i \leq 2} \det\nolimits_{\mathbf{Z}_\ell[G]}(1 - F_{*,i} \cdot u | \text{H}_i(X, \mathbf{Z}_\ell))^{(-1)^{i+1}}, \tag{36}$$

for any prime number $\ell \neq p$. Lemma 4.1.1 (1) and (2) and equality (36) give

$$\Theta(s) = \frac{\det\nolimits_{\mathbf{Z}_\ell[G]}(1 - F_{*,1} \cdot u | \text{H}_1(X, \mathbf{Z}_\ell))}{(1 - \sigma^{-1}u) \cdot (1 - q\sigma^{-1}u)}. \tag{37}$$

On the other hand, since $K/k$ is unramified everywhere, we have

$$\Theta_{S,T}(s) = \prod_{v \in S}(1 - \sigma_v^{-1}Nv^{-s}) \cdot \prod_{v \in T}(1 - \sigma_v^{-1}Nv^{1-s}) \cdot \Theta(s).$$

This equality and (37) give (1) in Lemma 4.2.1.

(2) This follows directly from (1) and Theorem 1.7.4.1 (3), with $i = 1$ and $Q(u) = Q_1(u)$. $\square$

**Proposition 4.2.2.** *Under the assumptions above, we have the following:*
  (1) $\Theta_{S,T}(0) \in \mathrm{Fitt}_{\mathbf{Z}[G]}(A_{S,T}) \cap \mathbf{Z}[G]_{0,S}$.
  (2) *If $\ell$ is a prime number, such that $\gcd(\ell, g) = 1$, or $\gcd(\ell, d_w) = 1$, for some $w \in S_K$, then $\Theta_{S,T}(0) \in \mathrm{Fitt}_{\mathbf{Z}_\ell[G]}(A_{S,T} \otimes \mathbf{Z}_\ell) \cdot \mathbf{Z}_\ell[G]_{0,S}$.*

**Proof.** Let us fix a prime $v_0 \in S$ and let $S_0 = \{v_0\}$. Obviously the set of data $(K/k, S, T, 0)$ still satisfies hypotheses (H). (Recall that $K/k$ is everywhere unramified.) There is a surjective morphism of $\mathbf{Z}[G]$–modules $A_{S_0,T} \longrightarrow A_{S,T}$ (see (9)), which gives a $\mathbf{Z}_\ell[G]$–surjection $A_{S_0,T} \otimes \mathbf{Z}_\ell \longrightarrow A_{S,T} \otimes \mathbf{Z}_\ell$, for every $\ell$ as above. These imply that

$$\mathrm{Fitt}_{\mathbf{Z}[G]}(A_{S_0,T}) \subseteq \mathrm{Fitt}_{\mathbf{Z}[G]}(A_{S,T}), \quad \mathrm{Fitt}_{\mathbf{Z}_\ell[G]}(A_{S_0,T} \otimes \mathbf{Z}_\ell) \subseteq \mathrm{Fitt}_{\mathbf{Z}_\ell[G]}(A_{S,T} \otimes \mathbf{Z}_\ell)$$

(see §1.4(f)). These relations together with

$$\Theta_{S,T}(0) = \prod_{v \in S \setminus \{v_0\}} (1 - \sigma_v^{-1}) \cdot \Theta_{S_0,T}(0), \text{ and } \prod_{v \in S \setminus \{v_0\}} (1 - \sigma_v^{-1}) \cdot \mathbf{Z}[G]_{0,S_0} \subseteq \mathbf{Z}[G]_{0,S},$$

imply that it is enough to prove the proposition above in the case $|S| = 1$. From now on we will assume that $S = \{v_0\}$, and we will fix a prime $w_0 \in S_K$. We will need three lemmas.

**Lemma 4.2.3.** *Let us consider the $\mathbf{Z}[G]$–module $\mathbf{Z}/d_{w_0}\mathbf{Z}$ with trivial $G$–action. Then*

$$\frac{1 - \sigma^{-d_{w_0}}}{1 - \sigma^{-1}} \in \mathrm{Fitt}_{\mathbf{Z}[G]}(\mathbf{Z}/d_{w_0}\mathbf{Z}).$$

**Proof.** Since $G$ acts trivially on $\mathbf{Z}/d_{w_0}\mathbf{Z}$, $\frac{1-\sigma^{-d_{w_0}}}{1-\sigma^{-1}} \in \mathrm{Ann}_{\mathbf{Z}[G]}(\mathbf{Z}/d_{w_0}\mathbf{Z})$. On the other hand, $\mathbf{Z}/d_{w_0}\mathbf{Z}$ is a cyclic $\mathbf{Z}[G]$–module, and therefore §1.4(a) shows that

$$\mathrm{Fitt}_{\mathbf{Z}[G]}(\mathbf{Z}/d_{w_0}\mathbf{Z}) = \mathrm{Ann}_{\mathbf{Z}[G]}(\mathbf{Z}/d_{w_0}\mathbf{Z}).$$

□

**Lemma 4.2.4.** *For any $v \in T$, the following hold true:*
  (1) $(1 - (q\sigma^{-1})^{d_v}) \in \mathrm{Fitt}_{\mathbf{Z}[G]}(\bigoplus_{w|v} \mathbf{F}_{q^\nu}(w)^\times)$.
  (2) $\frac{(1-(q\sigma^{-1})^{d_v})}{1-q\sigma^{-1}} \in \mathrm{Fitt}_{\mathbf{Z}[G]}\left((\bigoplus_{w|v} \mathbf{F}_{q^\nu}(w)^\times) / \mathbf{F}_{q^\nu}^\times\right),$
  *where $\mathbf{F}_{q^\nu}^\times$ is embeded in $\bigoplus_{w|v} \mathbf{F}_{q^\nu}(w)^\times$ diagonally ( $x \longrightarrow \bigoplus_{w|v}(x \mod w)$ ).*

**Proof.** Let us fix a prime $w_0$ in $K$, such that $w_0|v$. There is an isomorphism of $\mathbf{Z}[G]$–modules

$$\bigoplus_{w|v} \mathbf{F}_{q^\nu}(w)^\times \xrightarrow{\sim} \mathbf{F}_{q^\nu}(w_0)^\times \bigotimes_{\mathbf{Z}[G_v]} \mathbf{Z}[G],$$

where $G_v$ is the decomposition group of $v$ relative to $K/k$. Since $\mathbf{F}_{q^\nu}(w_0)^\times$ is a cyclic $\mathbf{Z}[G_v]$–module, both $\bigoplus_{w|v} \mathbf{F}_{q^\nu}(w)^\times$ and $(\bigoplus_{w|v} \mathbf{F}_{q^\nu}(w)^\times)/\mathbf{F}_{q^\nu}^\times$ are $\mathbf{Z}[G]$–cyclic and therefore their Fitting ideals over $\mathbf{Z}[G]$ are equal to their $\mathbf{Z}[G]$–annihilators (see §1.4(a)).

(1) The equality $\sigma_v = \sigma^{d_v}$ implies that

$$1 - (q\sigma^{-1})^{d_v} \in \mathrm{Ann}_{\mathbf{Z}[G_v]}(\mathbf{F}_{q^\nu}(w)^\times) \subseteq \mathrm{Ann}_{\mathbf{Z}[G]}(\bigoplus_{w|v} \mathbf{F}_{q^\nu}(w)^\times).$$

(2) Let $I = \{\nu, \nu-1, \ldots, \nu - r_v + 1\}$. Then $\{\sigma^i \mid i \in I\}$ is a set of coset representatives for the quotient $G/G_v$ and $\{w_0^{\sigma^i} \mid i \in I\} = \{w \mid w \text{ in } K, w|v\}$. Let

$$j_w : \mathbf{F}_{q^\nu}^\times \longrightarrow \mathbf{F}_{q^\nu}(w)^\times,$$

defined by $j_w(x) = x \bmod w$, for every $w|v$.

One can easily show that

$$j_{w_0^{\sigma^i}}(x) = \sigma^i[j_{w_0}(x^{q^{\nu-i}})], \forall\, i \in I,\, \forall x \in \mathbf{F}_{q^\nu}^\times. \tag{38}$$

We will identify $\mathbf{F}_{q^\nu}^\times$ with a subgroup of $\mathbf{F}_{q^\nu}(w_0)^\times$, via $j_{w_0}$. If $\zeta$ is a generator of $\mathbf{F}_{q^\nu}(w_0)^\times$ then, under this identification,

$$\mathbf{F}_{q^\nu}^\times = \langle \zeta^\alpha \rangle, \text{ where } \alpha = \frac{q^{\nu d_v/r_v} - 1}{q^\nu - 1}. \tag{39}$$

Relations (38) and (39) imply that the image of the composition of maps

$$J : \mathbf{F}_{q^\nu}^\times \xrightarrow{\oplus j_w}_{w|v} \bigoplus_{w|v} \mathbf{F}_{q^\nu}(w)^\times \xrightarrow{\sim} \mathbf{F}_{q^\nu}(w_0)^\times \bigotimes_{\mathbf{Z}[G_v]} \mathbf{Z}[G] \xrightarrow{\sim} \bigoplus_{i \in I} \sigma^i \cdot \mathbf{F}_{q^\nu}(w_0)^\times$$

satisfies the equality

$$\mathrm{Im}(J) = \{\sum_{i \in I} \sigma^i \cdot \zeta^{\alpha \cdot q^{\nu-i} \cdot n} \mid n \in \mathbf{Z}/(q^\nu - 1)\mathbf{Z}\}. \tag{40}$$

Let us also observe that $\bigoplus_{i \in I} \sigma^i \cdot \mathbf{F}_{q^\nu}(w_0)^\times$ is generated over $\mathbf{Z}[G]$ by the element $\overline{\zeta}$ having $\zeta = \sigma^\nu \zeta$ in the component corresponding to $i = \nu$, and $0$ in all the other components.

Let us consider the "$q^{d_v}$–evaluation map":

$$\phi_q : \mathbf{Z}[G_v] \longrightarrow \mathbf{Z}/(q^{d_v \nu/r_v} - 1)\mathbf{Z},$$

defined by $\phi_q(\sum_{j\geq 0} a_j \sigma_v{}^j) = \sum_{j\geq 0} a_j (q^{d_v})^j$. The equality $\sigma^{d_v}(\zeta) = \zeta^{q^{d_v}}$ shows that if $\overline{\alpha} = \sum_{i\in I} \sigma^i \cdot \alpha_i \in \mathbf{Z}[G]$, with $\alpha_i \in \mathbf{Z}[G_v]$, $\forall i \in I$, then $\overline{\alpha}\cdot\overline{\zeta} = \sum_{i\in I} \sigma^i \cdot \zeta^{\phi_q(\alpha_i)}$. In light of these observations, equality (40) shows that:

$$\mathrm{Ann}_{\mathbf{Z}[G]}\left((\bigoplus_{w|v}\mathbf{F}_{q^\nu}(w)^\times)/\mathbf{F}_{q^\nu}^\times\right) = \mathrm{Ann}_{\mathbf{Z}[G]}(\overline{\zeta}) =$$
$$= \{\sum_{i\in I} \sigma^i \cdot \alpha_i |\, \alpha_i \in \mathbf{Z}[G_v], \exists n \in \mathbf{Z}, \phi_q(\alpha_i) \equiv n \cdot \alpha q^{\nu-i} \mod (q^{\nu d_v/r_v} - 1), \forall i \in I\}. \quad (41)$$

One can write
$$\frac{1 - (q\sigma^{-1})^{d_v}}{1 - q\sigma^{-1}} = \sum_{i\in I} \sigma^i \cdot [q^{\nu-i} \cdot A],$$

where $A = \frac{1-(q\sigma^{-1})^{d_v}}{1-(q\sigma^{-1})^{r_v}}$. According to (41), in order to prove statement (2) in our lemma, we would simply have to show that $\alpha | \phi_q(A)$ in $\mathbf{Z}/(q^{d_v \nu/r_v} - 1)\mathbf{Z}$. This is an easy consequence of the fact that $\gcd(d_v/r_v, \nu/r_v) = 1$ (see comments at the beginning of §4.1). □

**Lemma 4.2.5.** *Let $\ell$ be any prime number. Then*

$$\det_{\mathbf{Z}_\ell[G]}(1 - F_{*,1}|\, \mathrm{H}_1(X, \mathbf{Z}_\ell)) \in \mathrm{Fitt}_{\mathbf{Z}_\ell[G]}\left(\mathrm{Pic}^0(X_1) \otimes \mathbf{Z}_\ell\right).$$

**Proof.** Lemma 3.1.5 applied to $M = \mathrm{Pic}^0(X)$ and $f = F_{*,1}$ gives an exact sequence of $\mathbf{Z}_\ell[G]$–modules

$$0 \longrightarrow \mathrm{T}_\ell(\mathrm{Pic}^0(X)) \xrightarrow{1-F_{*,1}} \mathrm{T}_\ell(\mathrm{Pic}^0(X)) \longrightarrow \mathrm{Pic}^0(X_1) \otimes \mathbf{Z}_\ell \longrightarrow 0. \quad (42)$$

(take into account that $\mathrm{Pic}^0(X)^{F_{*,1}} = \mathrm{Pic}^0(X_1)$). Lemma 4.1.1 (3) shows that $\mathrm{T}_\ell(\mathrm{Pic}^0(X)) = \mathrm{H}_1(X, \mathbf{Z}_\ell)$ is a free $\mathbf{Z}_\ell[G]$–module, for every prime number $\ell$. The exact sequence (42) gives therefore a free $\mathbf{Z}_\ell[G]$–resolution of $\mathrm{Pic}^0(X_1) \otimes \mathbf{Z}_\ell$. Lemma 4.2.5 is now a direct consequence of the definition of the Fitting ideal. □

We are now ready to conclude the proof of Proposition 4.2.2.

(1) By definition, $\Theta_{S,T}(0) \in \mathbf{Z}[G]_{0,S}$, and therefore all we have to prove is that $\Theta_{S,T}(0) \in \mathrm{Fitt}_{\mathbf{Z}[G]}(A_{S,T})$ which, according to §1.4(c), is equivalent to

$$\Theta_{S,T}(0) \in \mathrm{Fitt}_{\mathbf{Z}_\ell[G]}(A_{S,T} \otimes \mathbf{Z}_\ell), \forall \ell \text{ prime}.$$

Let us fix a prime number $\ell$. We will refer again to commutative diagram (11). With the same notations, we have an exact sequence of $\mathbf{Z}[G]$–modules

$$\mathrm{Pic}^0(X_1) \xrightarrow{\psi_S} A_S \longrightarrow \mathrm{coker}(\psi_S) \longrightarrow 0, \quad (43)$$

and $\operatorname{coker}(\phi_S) \xrightarrow{\sim} \operatorname{coker}(\psi_S)$ as $\mathbf{Z}[G]$–modules. The only difference is that, in our situation, due to the fact that $S = \{v_0\}$, $d_w = d_{w_0}$, for any $w \in S_K$. We therefore have a $\mathbf{Z}[G]$–isomorphism $\operatorname{coker}(\phi_S) \xrightarrow{\sim} \mathbf{Z}/d_{w_0}\mathbf{Z}$, which gives a $\mathbf{Z}[G]$–isomorphism

$$\operatorname{coker}(\psi_S) \xrightarrow{\sim} \mathbf{Z}/d_{w_0}\mathbf{Z}, \tag{44}$$

with $G$ acting trivially on $\mathbf{Z}/d_{w_0}\mathbf{Z}$. Lemmas 4.2.3 and 4.2.5, in addition to exact sequence (43) and property (e), §1.4, imply that

$$\frac{1-\sigma^{-d_{w_0}}}{1-\sigma^{-1}} \cdot \det\nolimits_{\mathbf{Z}_\ell[G]}(1 - F_{*,1}|\, \mathrm{H}_1(X, \mathbf{Z}_\ell)) \in \operatorname{Fitt}_{\mathbf{Z}_\ell[G]}(A_S \otimes \mathbf{Z}_\ell). \tag{45}$$

Let us fix a prime $v_T \in T$ and let $T_0 = \{v_T\}$. We have the following exact sequences of $\mathbf{Z}[G]$–modules:

$$0 \longrightarrow (\bigoplus_{w|v_T} \mathbf{F}_{q^\nu}(w)^\times)/(U_S/U_{S,T_0}) \longrightarrow A_{S,T_0} \longrightarrow A_S \longrightarrow 0,$$

$$(\bigoplus_{w|v_T} \mathbf{F}_{q^\nu}(w)^\times)/\mathbf{F}_{q^\nu}^\times \longrightarrow (\bigoplus_{w|v_T} \mathbf{F}_{q^\nu}(w)^\times)/(U_S/U_{S,T_0}) \longrightarrow 0.$$

(The first exact sequence above is another way of writing (1), with $T = T_0$, and the second one is a consequence of the obvious injection $\mathbf{F}_{q^\nu}^\times \hookrightarrow U_S/U_{S,T_0}$.) These exact sequences combined with (45), Lemma 4.2.4 (2) and §1.4(e), imply that

$$\frac{1-\sigma^{-d_{w_0}}}{1-\sigma^{-1}} \cdot \frac{1-(q\sigma^{-1})^{d_{v_T}}}{1-q\sigma^{-1}} \cdot \det\nolimits_{\mathbf{Z}_\ell[G]}(1 - F_{*,1}|\, \mathrm{H}_1(X, \mathbf{Z}_\ell)) \in \operatorname{Fitt}_{\mathbf{Z}_\ell[G]}(A_{S,T_0} \otimes \mathbf{Z}_\ell). \tag{46}$$

And finally, in order to pass from $T_0$ to $T$, let us consider the following exact sequences of $\mathbf{Z}[G]$–modules

$$0 \longrightarrow (\bigoplus_{w \in (T\setminus T_0)_K} \mathbf{F}_{q^\nu}(w)^\times)/(U_{S,T\setminus T_0}/U_{S,T}) \longrightarrow A_{S,T} \longrightarrow A_{S,T_0} \longrightarrow 0,$$

$$(\bigoplus_{w \in (T\setminus T_0)_K} \mathbf{F}_{q^\nu}(w)^\times) \longrightarrow (\bigoplus_{w \in (T\setminus T_0)_K} \mathbf{F}_{q^\nu}(w)^\times)/(U_{S,T\setminus T_0}/U_{S,T}) \longrightarrow 0.$$

These exact sequences, combined with (46), Lemma 4.2.4 (1) and §1.4(e), imply

$$\frac{1-\sigma^{-d_{w_0}}}{1-\sigma^{-1}} \cdot \frac{\prod_{v \in T} 1-(q\sigma^{-1})^{d_v}}{1-q\sigma^{-1}} \cdot \det\nolimits_{\mathbf{Z}_\ell[G]}(1 - F_{*,1}|\, \mathrm{H}_1(X, \mathbf{Z}_\ell)) \in \operatorname{Fitt}_{\mathbf{Z}_\ell[G]}(A_{S,T} \otimes \mathbf{Z}_\ell). \tag{47}$$

Let us denote by $\mathcal{L}$ the left hand–side of (47). Lemma 4.2.1, with $s = 0$ (and therefore $u = 1$), implies that

$$\Theta_{S,T}(0) = \begin{cases} \frac{1-\sigma^{-d_v}}{1-\sigma^{-d_{w_0}}} \cdot \mathcal{L}, & \text{if } \ell \neq p, \\ Q(1) \cdot \frac{1-\sigma^{-d_v}}{1-\sigma^{-d_{w_0}}} \cdot \mathcal{L}, & \text{if } \ell = p. \end{cases} \tag{48}$$

Since $d_{w_0} = d_v/r_v$ and $Q(1) \in \mathbf{Z}_p[G]$, (47) together with (48) show that

$$\Theta_{S,T}(0) \in \text{Fitt}_{\mathbf{Z}_\ell[G]}(A_{S,T} \otimes \mathbf{Z}_\ell),$$

for every prime number $\ell$. This, as we remarked earlier, concludes the proof of Proposition 4.2.2 (1).

(2) Let $\ell$ be a prime number as in Proposition 4.2.2 (2). Let us first notice that, if $\gcd(\ell, g) = 1$, then (2) is a direct consequence of (1). Indeed, in this situation

$$\mathbf{Z}_\ell[G] = \bigoplus_{\psi \in \widehat{G}(\mathbf{Q}_\ell)} D_\psi, \quad \mathbf{Z}_\ell[G]_{0,S} = \bigoplus_{\psi \in \widehat{G}(\mathbf{Q}_\ell, 0)} D_\psi,$$

which, together with $\text{Fitt}_{\mathbf{Z}_\ell[G]}(A_{S,T} \otimes \mathbf{Z}_\ell) = \oplus_\psi \text{Fitt}_{D_\psi}(A_{S,T} \otimes D_\psi)$, imply that

$$\text{Fitt}_{\mathbf{Z}_\ell[G]}(A_{S,T} \otimes \mathbf{Z}_\ell) \cdot \mathbf{Z}_\ell[G]_{0,S} = \text{Fitt}_{\mathbf{Z}_\ell[G]}(A_{S,T} \otimes \mathbf{Z}_\ell) \cap \mathbf{Z}_\ell[G]_{0,S}.$$

Let us then suppose that $\ell | g$ and $\gcd(\ell, d_{w_0}) = 1$. With notations as in the proof of (1), this implies that

$$\text{coker}(\psi_S) \otimes \mathbf{Z}_\ell = 0. \tag{49}$$

The exact sequence (43) and Lemma 4.2.5 imply that (45) can be rewritten as

$$\det_{\mathbf{Z}_\ell[G]}(1 - F_{*,1} | H_1(X, \mathbf{Z}_\ell)) \in \text{Fitt}_{\mathbf{Z}_\ell[G]}(A_S \otimes \mathbf{Z}_\ell),$$

which shows that in (45–47) above, the factor $\frac{1-\sigma^{-d_{w_0}}}{1-\sigma^{-1}}$ is not needed anymore. (47) and (48) therefore imply that

$$\Theta_{S,T}(0) \in \frac{1-\sigma^{-d_v}}{1-\sigma^{-1}} \cdot \text{Fitt}_{\mathbf{Z}_\ell[G]}(A_{S,T} \otimes \mathbf{Z}_\ell). \tag{50}$$

On the other hand, $\frac{1-\sigma^{-d_v}}{1-\sigma^{-1}} \in \mathbf{Z}[G]_{0,S}$. Indeed, if $\chi \in \widehat{G}(\mathbf{C}) \setminus \widehat{G}(\mathbf{C}, 0)$, then $\chi \neq 1_G$ and $\chi(\sigma_v) = 1$ (see (3) and take into account that $|S| = 1$.) This implies that $\chi(\frac{1-\sigma^{-d_v}}{1-\sigma^{-1}}) = 0$, which proves our assertion. This fact together with (50) imply Proposition 4.2.2 (2). □

Our next goal is the proof of the following:

**Proposition 4.2.6.** *Let $\ell$ be a prime number, such that $\ell \,|\, g$, and $\ell \,|\, d_w$, for every $w \in S_K$. Then*
$$\Theta_{S,T}(0) \in \text{Fitt}_{\mathbf{Z}_\ell[G]}(A_{S,T} \otimes \mathbf{Z}_\ell) \cdot \mathbf{Z}_\ell\,[G]_{0,S}\,.$$

An argument similar to the one given at the begining of the proof of Proposition 4.2.2 shows that it is enough to prove this statement in the case $|S| = 1$. We therefore assume that $S = \{v_0\}$, and fix $w_0 \in S_K$. In what follows, $\mathbf{C}_\ell$ denotes a completion of the algebraic closure of $\mathbf{Q}_\ell$ with respect to the $\ell$–adic valuation. $W$ denotes a finite ring-extension of $\mathbf{Z}_\ell$, inside $\mathbf{C}_\ell$, containing the values of all $\chi \in \widehat{G}(\mathbf{C}_\ell)$. We will extend scalars to $W$ and prove Proposition 4.2.6 inside $W[G]$.

Let $\ell^m$ be the exact power of $\ell$ dividing $g = |G|$. Let $L$ be the $\ell$–Sylow subgroup of $G$. There is a decomposition $G = \Delta \times L$, where $\Delta$ is the maximal subgroup of $G$, satisfying $\gcd(\ell, |\Delta|) = 1$. Let $K_L$ and $K_\Delta$ be the maximal subfields of $K$ fixed by $L$ and $\Delta$ respectively. Let $\overline{L} = G(K_\Delta/k)$ and $\overline{\Delta} = G(K_L/k)$. Restriction maps give isomorphisms $\Delta \xrightarrow{\sim} \overline{\Delta}$ and $L \xrightarrow{\sim} \overline{L}$.

**Remark 1.** Due to $d_{w_0} = d_v/\gcd(d_v, g)$, the fact that $\ell \,|\, d_{w_0}$ implies that $\ell^m \,|\, d_v$, and therefore $v$ splits completely in $K_\Delta$ (see §4.1). This implies that $\sigma_v \in \Delta$. Let us consider the decomposition $\widehat{G}(\mathbf{C}_\ell) = \widehat{\Delta}(\mathbf{C}_\ell) \times \widehat{L}(\mathbf{C}_\ell)$ and let $\chi$, $\delta$ and $\lambda$ denote generic elements of $\widehat{G}(\mathbf{C}_\ell)$, $\widehat{\Delta}(\mathbf{C}_\ell)$ and $\widehat{L}(\mathbf{C}_\ell)$ respectively. We can therefore conclude that a character $\chi = (\delta, \lambda)$ belongs to $\widehat{G}(\mathbf{C}_\ell, 0)$ if and only if $\delta(\sigma_v) \neq 1$, or $\chi = 1_G$ (i.e. $\delta = 1_\Delta$ and $\lambda = 1_L$) (see (3)).

Since $\gcd(\ell, |\Delta|) = 1$, we can decompose the group ring $W[G]$ into its $\delta$–components:
$$W[G] = W[\Delta][L] \xrightarrow{\underset{\delta \in \widehat{\Delta}(\mathbf{C}_\ell)}{\oplus \delta}} \bigoplus_{\delta \in \widehat{\Delta}(\mathbf{C}_\ell)} W_\delta[L],$$

where $W_\delta = W$ and $\delta(\sum_{\sigma \in \Delta, \tau \in L} a_{\sigma,\tau}\sigma\tau) = \sum_{\tau \in L}(\sum_{\sigma \in \Delta} a_{\sigma,\tau}\delta(\sigma)) \cdot \tau$, for every $\delta \in \widehat{\Delta}(\mathbf{C}_\ell)$. According to Remark 1 above, this induces an isomorphism
$$W[G]_{0,S} \xrightarrow{\sim} \left(\bigoplus_{\delta(\sigma_v) \neq 1} W_\delta[L]\right) \bigoplus W_{1_\Delta} \cdot N_L,$$

where $N_L = \sum_{\sigma \in L} \sigma$. We thus have a $W[G]$–isomorphism
$$W[G]_{0,S} \cdot \text{Fitt}_{W[G]}(A_{S,T} \otimes W) \xrightarrow{\sim} \left(\bigoplus_{\delta(\sigma_v) \neq 1} \text{Fitt}_{W_\delta[L]}(A_{S,T} \otimes W)^\delta\right) \oplus$$
$$\oplus \text{Fitt}_{W[L]}(A_{S,T} \otimes W)^\Delta \cdot N_L.$$

This implies that, in order to prove
$$\Theta_{S,T}(0) \in W[G]_{0,S} \cdot \text{Fitt}_{W[G]}(A_{S,T} \otimes W),$$

it is necessary and sufficient to show the following:
  (i)   If $\delta(\sigma_v) = 1$, $\delta \ne 1_\Delta$, then $\delta(\Theta_{S,T}(0)) = 0$.
  (ii)  If $\delta(\sigma_v) \ne 1$, then $\delta(\Theta_{S,T}(0)) \in \operatorname{Fitt}_{W_\delta[L]}(A_{S,T} \otimes W)^\delta$.
  (iii) $1_\Delta(\Theta_{S,T}(0)) \in \operatorname{Fitt}_{W[L]}(A_{S,T} \otimes W)^\Delta \cdot N_L$.

Assertions (i) and (ii) above are direct consequences of Proposition 4.2.2 (1), after extending scalars to $W$. Therefore, in order to prove Proposition 4.2.6, one has to prove (iii) above.

**Proof of Prop. 4.2.6.** In the following arguments, if $F$ is an intermediate field of $K/k$, then $A_{F,S,T}$ denotes its $(S,T)$–ideal class group, as defined in §1.1.

The fact that $\Theta_{S,T}(s) \in W[G]_{0,S}$, combined with the class-number formula (2), implies that

$$1_\Delta(\Theta_{S,T}(0)) = \zeta_{k,S,T}(0) \cdot \frac{N_L}{|L|} = \frac{|A_{k,S,T}|}{|L|} \cdot N_L.$$

Via the isomorphism $W[L] \xrightarrow{\sim} W[\overline{L}]$, we can therefore say that, inside $W[\overline{L}]$,

$$1_\Delta(\Theta_{S,T}(0)) = \frac{|A_{k,S,T}|}{|\overline{L}|} \cdot N_{\overline{L}}, \tag{51}$$

where $N_{\overline{L}} = \sum_{\sigma \in \overline{L}} \sigma$. The fact that $\gcd(\ell, |\Delta|) = 1$, combined with a standard class-field theoretical argument based on the class-field interpretation of the groups $A_{F,S,T}$ (see §1.1), implies that we have an isomorphism of $W[\overline{L}]$–modules

$$(A_{S,T} \otimes W)^\Delta \xrightarrow{\sim} A_{K_\Delta,S,T} \otimes W.$$

This obviously implies that

$$\operatorname{Fitt}_{W[\overline{L}]}(A_{S,T} \otimes W)^\Delta = \operatorname{Fitt}_{W[\overline{L}]}(A_{K_\Delta,S,T} \otimes W). \tag{52}$$

Equalities (51) and (52) above show that (iii) would be a consequence (after extending scalars to $W$) of the following

**Lemma 4.2.7.** *With notations as above,*

$$\frac{|A_{k,S,T}|}{|\overline{L}|} \cdot N_{\overline{L}} \in \operatorname{Fitt}_{\mathbf{Z}[\overline{L}]}(A_{K_\Delta,S,T}) \cdot N_{\overline{L}}.$$

**Proof of Lemma 4.2.7.** Let us consider the exact sequence

$$0 \longrightarrow I_{\overline{L}} \longrightarrow \mathbf{Z}[\overline{L}] \xrightarrow{s} \mathbf{Z} \longrightarrow 0,$$

where $s(\sum_{\sigma \in \overline{L}} a_\sigma \cdot \sigma) = \sum_{\sigma \in \overline{L}} a_\sigma$, for every $\sum_{\sigma \in \overline{L}} a_\sigma \cdot \sigma \in \mathbf{Z}[\overline{L}]$, and $I_{\overline{L}} = \ker(s)$.

Obviously $I_{\overline{L}} = \mathrm{Ann}_{\mathbf{Z}[\overline{L}]}(N_{\overline{L}})$, and therefore the statement in Lemma 4.2.7 is equivalent to

$$\frac{|A_{k,S,T}|}{|\overline{L}|} \in s(\mathrm{Fitt}_{\mathbf{Z}[\overline{L}]}(A_{K_\Delta,S,T})).$$

On the other hand, properties (b) and (f), §1.4 imply that

$$s(\mathrm{Fitt}_{\mathbf{Z}[\overline{L}]}(A_{K_\Delta,S,T})) = \mathrm{Fitt}_{\mathbf{Z}}(A_{K_\Delta,S,T}/I_{\overline{L}}A_{K_\Delta,S,T}) = |A_{K_\Delta,S,T}/I_{\overline{L}}A_{K_\Delta,S,T}| \cdot \mathbf{Z}.$$

The statement in Lemma 4.2.7 is therefore equivalent to

$$|A_{K_\Delta,S,T}/I_{\overline{L}}A_{K_\Delta,S,T}| \mid \frac{|A_{k,S,T}|}{|\overline{L}|}.$$

The next lemma shows that something even stronger holds true.

**Lemma 4.2.8.** *With notations as above,*

$$|A_{K_\Delta,S,T}/I_{\overline{L}}A_{K_\Delta,S,T}| = \frac{|A_{k,S,T}|}{|\overline{L}|}.$$

**Proof of Lemma 4.2.8.** Let $K_{\Delta,S,T}$ and $k_{S,T}$ be the $(S,T)$–ray class-fields of $K_\Delta$ and $k$ respectively. We remind the reader that these are the maximal abelian extensions of $K_\Delta$ and $k$ respectively, which are tamely ramified at primes above $T$ and unramified outside $T$, and in which all primes in $S$ split completely (see class-field theoretical interpretation of $A_{S,T}$ in §1.1). The corresponding Artin reciprocity maps give the following group isomorphisms:

$$A_{K_\Delta,S,T} \xrightarrow{\sim} G(K_{\Delta,S,T}/K_\Delta) \tag{53}$$

$$A_{k,S,T} \xrightarrow{\sim} G(k_{S,T}/k). \tag{54}$$

The maximality of $K_{\Delta,S,T}$, the fact that $v_0$ splits completely in $K_\Delta$ (see Remark 1), and that $K_\Delta/k$ is everywhere unramified, imply that $K_\Delta \subseteq k_{S,T} \subseteq K_{\Delta,S,T}$ and that $K_{\Delta,S,T}/k$ is Galois. Let $\mathcal{G} = G(K_{\Delta,S,T}/k)$. We have an exact sequence of groups

$$1 \longrightarrow G(K_{\Delta,S,T}/K_\Delta) \longrightarrow \mathcal{G} \longrightarrow \overline{L} \longrightarrow 1. \tag{55}$$

Since $G(K_{\Delta,S,T}/K_\Delta)$ is abelian, $\overline{L}$ acts on it by (lift and) conjugation and this way (53) becomes an isomorphism of $\mathbf{Z}[\overline{L}]$–modules.

Let $[\mathcal{G},\mathcal{G}]$ denote the group theoretical commutator of $\mathcal{G}$. The maximality of $k_{S,T}$ implies that

$$[\mathcal{G},\mathcal{G}] = G(K_{\Delta,S,T}/k_{S,T}). \tag{56}$$

An easy group theoretical argument, based on the exact sequence (55) and on the fact that $\overline{L}$ is cyclic, shows that $[\mathcal{G},\mathcal{G}]$ is generated by elements of the form $\gamma^{\sigma-1}$, with $\gamma \in G(K_{\Delta,S,T}/K_\Delta)$ and $\sigma \in \overline{L}$. On the other hand, the ideal $I_{\overline{L}}$ is generated over $\mathbf{Z}$ by elements of the form $\sigma - 1$. This observation and (56) imply that, under the Artin isomorphism (53), we have

$$I_{\overline{L}} \cdot A_{K_\Delta,S,T} \xrightarrow{\sim} G(K_{\Delta,S,T}/k_{S,T}).$$

We therefore have the equalities

$$|A_{K_\Delta,S,T}/I_{\overline{L}} A_{K_\Delta,S,T}| = \frac{|G(K_{\Delta,S,T}/K_\Delta)|}{|G(K_{\Delta,S,T}/k_{S,T})|} = \frac{|G(k_{S,T}/k)|}{|\overline{L}|} = \frac{|A_{k,S,T}|}{|\overline{L}|}.$$

This concludes the proof of Lemma 4.2.8, Lemma 4.2.7 and Proposition 4.2.6. □

**Theorem 4.2.9.** *Let $K/k$ be a finite constant field extension of function fields, and let us suppose that the set of data $(K/k, S, T, 0)$ satisfies hypotheses (H). Then*

$$\Theta_{S,T}(0) \in \mathrm{Fitt}_{\mathbf{Z}[G]}(A_{S,T}) \cdot \mathbf{Z}\,[G]_{0,S}.$$

**Proof.** This statement is obviously equivalent to

$$\Theta_{S,T}(0) \in \mathrm{Fitt}_{\mathbf{Z}_\ell[G]}(A_{S,T} \otimes \mathbf{Z}_\ell) \cdot \mathbf{Z}_\ell\,[G]_{0,S}, \forall\, \ell \text{ prime,}$$

which follows from Propositions 4.2.2(2) and 4.2.6 above. □

4.3. ARBITRARY $r$

**Theorem 4.3.1.** *Let $K/k$ be a finite constant field extension of function fields, and let us suppose that the set of data $(K/k, S, T, r)$ satisfies hypotheses (H). Then there exists a unique element*

$$\varepsilon_{S,T} \in \mathbf{Z}\,[G]_{r,S} \cdot (\overset{r}{\wedge} U_{S,T}),$$

*such that $\Theta_{S,T}(0) = R_W(\varepsilon_{S,T})$.*

**Proof.** Let $\{v_1, \cdots, v_r\}$ be the distinguished set of $r$ distinct primes in $S$, which split completely in $K/k$. Let $S_0 = S \setminus \{v_1, \cdots, v_r\}$. The set of data $(K/k, S_0, T, 0)$ satisfies hypotheses (H) as well. Theorem 4.2.9 therefore implies that

$$\Theta_{S_0,T}(0) \in \mathrm{Fitt}_{\mathbf{Z}[G]}(A_{S_0,T}) \cdot \mathbf{Z}\,[G]_{0,S_0}.$$

The statement in Theorem 4.3.1 is now a direct consequence of the relation above, and Corollary 2.2 (2), with $S' = S$ and $S = S_0$. □

**Remark.** The inclusion

$$\mathbf{Z}\,[G]_{r,S} \cdot (\overset{r}{\wedge} U_{S,T}) \subseteq \Lambda_{S,T},$$

shows that Theorem 4.2.9 settles Conjecture B, for arbitrary $r$, under the assumption that $K/k$ is a constant field extension. The inclusion above is strict in most cases, and therefore Theorem 4.2.9 provides a refinement of Conjecture B, under the present assumptions.

**Appendix**

The goal of this appendix is the proof of Theorem 1.7.4.1. With notations as in §1.7, let us recall that, for any prime number $\ell$, and any $i = 0, 1, 2$,

$$P_{i,\ell}(u) \stackrel{\text{def}}{=} \det_{\mathbf{Q}_\ell[G]} \left(1 - F_{*,i} \cdot u | \, \mathrm{H}_i\left(X, \mathbf{Q}_\ell\right)\right),$$

where $\mathrm{H}_i(X, \mathbf{Q}_\ell)$ is the $i$–th $\ell$–adic homology group with coefficients in $\mathbf{Q}_\ell$, defined in §1.7. The statement we would like to prove is the following:

**Theorem.** *For any $i = 0, 1, 2$:*
  (1) *$P_{i,p}(u) \in \mathbf{Z}_p\left[1/g\right][G][u]$, and if $\ell \neq p$, then $P_{i,\ell} \in \mathbf{Z}\left[1/g\right][G][u]$.*
  (2) *If $\ell \neq p$, then $P_{i,\ell}(u)$ does not depend on $\ell$.*
     *Let $P_i(u) \stackrel{\text{def}}{=} P_{i,\ell}(u)$, for any $\ell \neq p$. Then:*
  (3) *There exist polynomials $Q_i(u) \in \mathbf{Z}_p\left[1/g\right][G][u]$ such that*

$$P_i(u) = P_{i,p}(u) \cdot Q_i(u).$$

  (4) *If $K/k$ is a constant field extension, then $Q_i(u) \in \mathbf{Z}_p[G][u]$.*

Before proceeding to the proof of the Theorem above, we need to summarize a few considerations made in [17], Chpt.V, on the Stickelberger function $\Theta(s)$.

Let $|Y_1|$ be the set of closed points of the $\mathbf{F}_q$–scheme $Y_1$. These are in one–to–one correspondence with the primes $v$ of the base field $k$. For every $v \in |Y_1|$, let $I_v$ be its innertia group in $G$, let $\overline{\sigma_v} \in G_v/I_v$ be its Frobenius class, and let

$$F_v \stackrel{\text{def}}{=} \frac{1}{|I_v|} \sum_{\tau \in \overline{\sigma_v}} \tau^{-1} \in \mathbf{Z}\left[1/g\right][G].$$

Then, if $u = q^{-s}$, we have the following equality

$$\Theta(u) = \prod_{v \in |Y_1|} \left(1 - F_v \cdot u^{-\deg v}\right) \in 1 + u\mathbf{Z}\left[1/g\right][G][[u]], \qquad (\text{a.1})$$

valid for $|u| < q^{-1}$ (which corresponds to $\mathrm{Re}(s) > 1$).

**Remark 1.** Let us observe that, if $K/k$ is a constant field extension, then $|I_v| = 1$, and therefore $F_v \in \mathbf{Z}[G]$, for every $v \in |Y_1|$. Equality (a.1) therefore shows that, under the present hypothesis,

$$\Theta(u) \in \mathbf{Z}[G][[u]]. \qquad (\text{a.2})$$

For every $v \in |Y_1|$, let $X_v$ be the finite set of closed points of the $\mathbf{F}$–scheme $X$ above $v$, and let $\mathrm{H}_0(X_v, \mathbf{Q})$ be the $\mathbf{Q}$–vector space generated by $X_v$, with the obvious $G$–action and $q$–power geometric Frobenius morphism action $F_*$. Lemma 2.1 in [17], Chpt.V, shows that

$$1 - F_v u^{-\deg v} = \det_{\mathbf{Q}[G]}\left(1 - F_* \cdot u | \operatorname{H}_0\left(X_v, \mathbf{Q}\right)\right), \quad \forall v \in |Y_1|.$$

If we combine these equalities with (a.1), we obtain

$$\Theta(u) = \prod_{v \in |Y_1|} \det_{\mathbf{Q}[G]}\left(1 - F_* \cdot u | \operatorname{H}_0\left(X_v, \mathbf{Q}\right)\right). \tag{a.3}$$

We are now prepared to prove (1) and (2) in the Theorem above.

**Proof of (1).** Since $\operatorname{H}_i(X, \mathbf{Z}_p)$ is a free $\mathbf{Z}_p$–module,

$$\operatorname{H}_i(X, \mathbf{Z}_p[1/g]) = \operatorname{H}_i(X, \mathbf{Z}_p) \otimes \mathbf{Z}_p[1/g]$$

is a projective $\mathbf{Z}_p[1/g][G]$–module (see Remark 1, §1.3), and therefore

$$P_{i,p}(u) = \det_{\mathbf{Z}_p[1/g][G]}\left(1 - F_{*,i} | \operatorname{H}_i(X, \mathbf{Z}_p[1/g])\right) \in \mathbf{Z}_p[1/g][G][u].$$

Let $\ell \neq p$ be a prime number. Since $\operatorname{H}_0(X, \mathbf{Z}[1/g])$ is $\mathbf{Z}[1/g][G]$–projective and $\operatorname{H}_0(X, \mathbf{Q}_\ell) = \operatorname{H}_0(X, \mathbf{Z}[1/g]) \bigotimes_{\mathbf{Z}[1/g]} \mathbf{Q}_\ell$, we have the equality

$$P_{0,\ell}(u) = \det_{\mathbf{Z}[1/g][G]}\left(1 - F_{*,0} \cdot u | \operatorname{H}_0(X, \mathbf{Z}[1/g])\right). \tag{a.4}$$

This obviously implies that $P_{0,\ell} \in \mathbf{Z}[1/g][G][u]$.

The $\mathbf{Q}_\ell[G]$–isomorphism $\operatorname{H}_2(X, \mathbf{Q}_\ell) \xrightarrow{\sim} \operatorname{H}_0(X, \mathbf{Q}_\ell)$, carying the action of $F_{*,2}$ into the action of $qF_{*,0}$ (see §1.7.3), gives the equality

$$P_{2,\ell}(u) = \det_{\mathbf{Q}_\ell[G]}\left(1 - qF_{*,0} \cdot u | \operatorname{H}_0(X, \mathbf{Q}_\ell)\right), \tag{a.5}$$

which, according to the remarks above, implies that

$$P_{2,\ell}(u) = \det_{\mathbf{Z}[1/g][G]}\left(1 - qF_{*,0} \cdot u | \operatorname{H}_0(X, \mathbf{Z}[1/g])\right). \tag{a.6}$$

The last equality obviously shows that $P_{2,\ell}(u) \in \mathbf{Z}[1/g][G][u]$.

Theorem 1.7.4.2 shows that

$$P_{1,\ell}(u) = \Theta(u) \cdot P_{0,\ell}(u) \cdot P_{2,\ell}(u) \tag{a.7}$$

and since $\Theta(u) \in \mathbf{Z}[1/g][G][[u]]$ (see (a.1)), we obviously have $P_{1,\ell} \in \mathbf{Z}[1/g][G][u]$ as well. $\square$

**Proof of (2).** Equalities (a.4) and (a.6) obviously show that $P_{0,\ell}(u)$ and $P_{1,\ell}(u)$ do not depend on $\ell$, $\ell \neq p$. Equality (a.7) therefore shows that $P_{1,\ell}(u)$ is independent of $\ell$, $\ell \neq p$, as well. $\square$

Before proceeding to the proof of statements (3) and (4) in the Theorem above, we will need to make a few considerations on the $p$–adic étale and the crystalline cohomology groups of $X$.

Let $\mathcal{H}(X) = \bigoplus_{i=0,1,2} \mathcal{H}^i(X)$ be either the $\ell$–adic étale cohomology of $X$ with coefficients in $\mathbf{C}_\ell$, $\bigoplus_{i=0,1,2} \mathrm{H}^i(X, \mathbf{C}_\ell)$, $\ell \neq p$, or the crystalline cohomology of $X$ with coefficients in $\mathbf{C}_p$, $\bigoplus_{i=0,1,2} \mathrm{H}^i_{\mathrm{cris}}(X, \mathbf{C}_p)$. Here $\mathbf{C}_\ell$ denotes the completion of the algebraic closure of $\mathbf{Q}_\ell$ with respect to the extension of the normalized $\ell$–adic valuation, for any prime number $\ell$. The definition of $\mathrm{H}^i(X, \mathbf{C}_\ell)$ is

$$\mathrm{H}^i(X, \mathbf{C}_\ell) = \mathrm{H}^i(X, \mathbf{Z}_\ell) \bigotimes_{\mathbf{Z}_\ell} \mathbf{C}_\ell,$$

where $\mathrm{H}^i(X, \mathbf{Z}_\ell) = \mathrm{H}_i(X, \mathbf{Z}_\ell)^*$ (functorial dual), for all $i = 0, 1, 2$, and every $\ell \neq p$ (see [10]).

Let $W(\mathbf{F})$ be the Witt vector ring of $\mathbf{F}$. Then

$$\mathrm{H}^i_{\mathrm{cris}}(X, \mathbf{C}_p) \stackrel{\mathrm{def}}{=} \mathrm{H}^i_{\mathrm{cris}}(X/W(\mathbf{F})) \bigotimes_{W(\mathbf{F})} \mathbf{C}_p,$$

where $\mathrm{H}^i_{\mathrm{cris}}(X/W(\mathbf{F}))$ are finitely generated $W(\mathbf{F})$–modules, functorial in the pair $X/\mathbf{F}$. For a precise definition of the modules $\mathrm{H}^i_{\mathrm{cris}}(X/W(\mathbf{F}))$ the reader can consult [1] and [2]. We will not need it for the present considerations.

In what follows $\mathbb{K}$ will denote either $\mathbf{C}_p$ or $\mathbf{C}_\ell$, $(\ell \neq p)$, depending on which cohomology theory is being used. The correspondences $X \longrightarrow \mathcal{H}^i(X)$ are contravariant functors from the category of smooth, projective $\mathbf{F}$–schemes, to the category of $\mathbb{K}$–vector spaces. For every morphism of $\mathbf{F}$–schemes $f: X \longrightarrow X'$, we will denote by $f^*$ the $\mathbb{K}$–linear maps

$$\mathcal{H}^i(X') \xrightarrow{\mathcal{H}^i(f)} \mathcal{H}^i(X),$$

induced by $f$ at the level of cohomology. In particular, every $\sigma \in G$ induces an isomorphism of $\mathbf{F}$–schemes $\sigma_X: X \longrightarrow X$, and therefore isomorphisms of $\mathbb{K}$–vector spaces

$$\sigma_X^*: \mathcal{H}^i(X) \longrightarrow \mathcal{H}^i(X).$$

We can define an action of $G$ on $\mathcal{H}(X)$ by setting

$$\sigma \cdot h \stackrel{\mathrm{def}}{=} (\sigma_X^*)^{-1}(h),$$

for all $\sigma \in G$ and $h \in \mathcal{H}(X)$. This way the étale cohomology groups $\mathrm{H}^i(X, \mathbf{C}_\ell)$ are dual as Galois modules to the homology groups $\mathrm{H}_i(X, \mathbf{C}_\ell)$ defined in §1.7.4. In particular, for any character $\chi \in \widehat{G}(\mathbf{C}_\ell)$ we have $\mathbf{C}_\ell$–linear isomorphisms

$$\mathrm{H}^i(X, \mathbf{C}_\ell)^\chi \stackrel{\sim}{\longrightarrow} \mathrm{H}_i(X, \mathbf{C}_\ell)^{\chi^{-1}}, \tag{a.8}$$

carrying the action $F^*$ of the $q$–power geometric Frobenius morphism on the left hand–side into its action $F_{*,i}$ on the right hand–side.

Let $P^i_{\mathcal{H}}(u) \stackrel{\text{def}}{=} \det_{\mathbb{K}}\left(1 - F^* \cdot u \,|\, \mathcal{H}^i(X)\right)$. Theorem 1 in [8] shows that the polynomials $P^i_{\mathcal{H}}(u)$ do not depend on the cohomology theory $\mathcal{H}(X)$, and therefore we can make the notation $P^i(u) = P^i_{\mathcal{H}}(u)$, for all $i = 0, 1, 2$. The polynomials $P^i(u)$ have coefficients in $\mathbf{Z}$ and their reciprocal roots $\{\alpha_{ij}\}$ are algebraic integers, satisfying the relations

$$|\alpha_{ij}| = q^{i/2}, \qquad (a.9)$$

for all $i = 0, 1, 2$, and all $j = 1, \ldots, \dim_{\mathbb{K}} \mathcal{H}^i(X)$ (see [8]). Relation (a.9) is the dimension 1 case of Riemann's Hypothesis, proved by Deligne for smooth, projective $\mathbf{F}_q$–varieties of any dimension.

Let $P^{i,\chi}_{\mathcal{H}}(u) = \det_{\mathbb{K}}\left(1 - F^* \cdot u \,|\, \mathcal{H}^i(X)^\chi\right)$, for all $i = 0, 1, 2$, and $\chi \in \widehat{G}(\mathbb{K})$. We obviously have the relations

$$P^i(u) = \prod_{\chi \in \widehat{G}(\mathbb{K})} P^{i,\chi}_{\mathcal{H}}(u), \qquad (a.10)$$

for every $i = 0, 1, 2$, and every cohomology theory $\mathcal{H}(X)$ in our list.

**Remark 2.** Relations (a.9) and (a.10) imply that for a given cohomology theory $\mathcal{H}(X)$ and a given $\chi \in \widehat{G}(\mathbb{K})$, the polynomials $P^{i,\chi}_{\mathcal{H}}(u)$ have mutually disjoint sets of roots, for distinct values of $i$.

Our next goal is to prove the following:

**Proposition A.1.** *The polynomials $P^{i,\chi}_{\mathcal{H}}(u)$ do not depend on the cohomology theory $\mathcal{H}(X)$.*

Before proceeding to proving the statement above, we need an elementary lemma and a few remarks.

**Lemma A.2.** *If $V$ is a finitely generated $\mathbf{Q}[G]$–module, and $f \in \mathrm{End}_{\mathbf{Q}[G]}(V)$, then*

(1) $\mathrm{Tr}_{\mathbf{Q}[G]}(f; V) = 1/g \sum_{\sigma \in G} \mathrm{Tr}_{\mathbf{Q}}\left(\sigma^{-1} f; V\right) \cdot \sigma$.

(2) $\mathrm{Tr}_{\mathbb{K}}\left(f; (V \otimes \mathbb{K})^\chi\right) = 1/g \sum_{\sigma \in G} \mathrm{Tr}_{\mathbb{K}}\left(\sigma^{-1} f; V \otimes \mathbb{K}\right) \cdot \chi(\sigma)$, *for all $\chi \in \widehat{G}(\mathbb{K})$.*

**Proof.** See [17], Chpt. V, Lemme 2.6. □

In what follows we will freely use the mutually inverse isomorphisms

$$1 + u\mathbf{Q}[G][[u]] \xrightleftharpoons[\exp]{\log} u\mathbf{Q}[G][[u]],$$

and also the identity

$$\log \det_{\mathbf{Q}[G]}(1 - f \cdot u \,|\, V) = -\sum_{n \geq 1} \mathrm{Tr}_{\mathbf{Q}[G]}(f^n; V) \cdot u^n/n, \qquad (a.11)$$

for any $V$ and $f$ as in Lemma A.2 (see [17], Chpt.V).

Any cohomology theory $\mathcal{H}(X)$ in the list above satisfies the following form of the Lefschetz fixed point formula :

**Theorem A.3.** *Let $f : X \longrightarrow X$ be a morphism of $\mathbf{F}$–schemes, with isolated, multiplicity 1 fixed points. Let $\Lambda(f, X)$ be the number of fixed points of $f$. Then*

$$\Lambda(f, X) = \sum_{i=0,1,2} \mathrm{Tr}_{\mathbb{K}}\left(f^* \,;\, \mathcal{H}^i(X)\right) .$$

For a proof of the theorem above in the case of étale cohomology, the reader can consult [10]. The case of crystalline cohomology is treated in [1] and [2]. We now recall the definition of a fixed point of multiplicity 1 (see [7], Appendix C).

**Definition A.4.** A closed point $x \in X$, fixed by $f$, has multiplicity 1 if the map

$$\Omega_{X,x} \xrightarrow{1-df_x} \Omega_{X,x} ,$$

induced by $f$ at the level of the $x$–germ $\Omega_{X,x}$ of the sheaf of differentials $\Omega_X$ of $X$, is injective.

We are now prepared to prove Proposition A.1.

**Proof of Proposition A.1.** (Compare [17], Chpt.V, Thm.2.5) According to (a.3) and (a.11), we have

$$\log \Theta(u) = \sum_{v \in |Y_1|} \log \det\nolimits_{\mathbf{Q}[G]} \left(1 - F_* \cdot u |\, \mathrm{H}_0(X_v, \mathbf{Q})\right) =$$
$$= \sum_{v,n} \mathrm{Tr}_{\mathbf{Q}[G]}\left(F_*^n \,;\, \mathrm{H}_0(X_v, \mathbf{Q})\right) \cdot u^n/n .$$

Lemma A.2 (1) shows that

$$\mathrm{Tr}_{\mathbf{Q}[G]}\left(F_*^n \,;\, \mathrm{H}_0(X_v, \mathbf{Q})\right) = 1/g \mathrm{Tr}_{\mathbf{Q}}\left(\sigma_*^{-1} F_*^n \,;\, \mathrm{H}_0(X_v, \mathbf{Q})\right) .$$

But since $\sigma_*^{-1} F_*^n$ acts on $\mathrm{H}_0(X_v, \mathbf{Q})$ by simply permuting the points of $X_v$, we also have the relation

$$\mathrm{Tr}_{\mathbf{Q}}\left(\sigma_*^{-1} F_*^n \,;\, \mathrm{H}_0(X_v, \mathbf{Q})\right) = \Lambda\left(\sigma_X^{-1} F^n \,;\, X_v\right) ,$$

where the left hand–side in the last equality represents the number of fixed points of $\sigma_X^{-1} F^n|_{X_v}$. This implies that

$$\log \Theta(u) = 1/g \sum_{v,n,\sigma} \Lambda\left(\sigma_X^{-1} F^n \,;\, X_v\right) \cdot \sigma \cdot u^n/n ,$$

and since $X = \coprod_{v \in |Y_1|} X_v$, we obtain

$$\log \Theta(u) = 1/g \sum_{n,\sigma} \Lambda\left(\sigma_X^{-1} F^n \,;\, X\right) \cdot \sigma \cdot u^n/n \,,$$

where $\Lambda\left(\sigma_X^{-1} F^n \,;\, X\right)$ is the number of fixed points of the **F**–scheme morphism $\sigma_X^{-1} F^n : X \longrightarrow X$, for all $n \geq 1$.

Let $x \in X$ be a point fixed by $\sigma_X^{-1} F^n$, for some $n \geq 1$. One can easily show that $x$ is isolated ($\dim X = 1$), and that it has multiplicity 1, in the sense of Definition A.4. We are therefore entitled to use Theorem A.3 in order to compute the numbers $\Lambda\left(\sigma_X^{-1} F^n \,;\, X\right)$, for all $n \geq 1$. We thus have

$$\log \Theta(u) = 1/g \sum_{n,\sigma,i} (-1)^i \text{Tr}_{\mathbb{K}} \left( \left(\sigma_X^{-1} F^n\right)^* \,;\, \mathcal{H}^i(X) \right) \cdot \sigma \cdot u^n/n$$

$$= 1/g \sum_{n,\sigma,i} (-1)^i \text{Tr}_{\mathbb{K}} \left( \sigma F^{*n} \,;\, \mathcal{H}^i(X) \right) \cdot \sigma \cdot u^n/n \,.$$

This implies that, for any $\chi \in \widehat{G}(\mathbb{K})$, we have

$$\log L(u, \chi) = \log\left(\chi^{-1} \Theta(u)\right) = \chi^{-1}\left(\log \Theta(u)\right) =$$
$$= 1/g \sum_{n,\sigma,i} (-1)^i \text{Tr}_{\mathbb{K}} \left( \sigma^{-1} F^{*n} \,;\, \mathcal{H}^i(X) \right) \cdot \chi(\sigma) u^n/n \,.$$

Let us now combine the last equality with (a.11) and with Lemma A.2 (2), for $V = \mathcal{H}^i(X)$ and $f = F^{*n}$. We obtain

$$\log L(u, \chi) = \sum_{n,i} (-1)^i \text{Tr}_{\mathbb{K}} \left( F^{*n} \,;\, \mathcal{H}^i(X)^\chi \right) \cdot u^n/n =$$
$$= \sum_i (-1)^{i+1} \log \det_{\mathbb{K}} \left( 1 - F^* \cdot u \,;\, \mathcal{H}^i(X)^\chi \right) =$$
$$= \log \left[ \prod_i \det_{\mathbb{K}} \left( 1 - F^* \cdot u \,;\, \mathcal{H}^i(X)^\chi \right)^{(-1)^{i+1}} \right] =$$
$$= \log \left[ \prod_i P_{\mathcal{H}}^{i,\chi}(u)^{(-1)^{i+1}} \right] \,.$$

And now, by taking *exp* of both sides of the last equality, we obtain:

$$L(u, \chi) = \prod_i P_{\mathcal{H}}^{i,\chi}(u)^{(-1)^{i+1}} \,,$$

for all $\chi$ and $\mathcal{H}(X)$ as above. This shows that the product $\prod_i P_{\mathcal{H}}^{i,\chi}(u)^{(-1)^{i+1}}$ does not depend on the cohomology theory $\mathcal{H}(X)$. Remark 2 therefore shows that the individual polynomials $P_{\mathcal{H}}^{i,\chi}(u)$ do not depend on $\mathcal{H}(X)$, which concludes the proof of Proposition A.1. □

Let us make the notation $P^{i,\chi}(u) = P_{\mathcal{H}}^{i,\chi}(u)$, for all $i = 0, 1, 2, \chi \in \widehat{G}$, and some cohomology theory $\mathcal{H}(X)$ in our list. The isomorphism (a.8) shows that, for all prime numbers $\ell \neq p$, we have

$$P^{i,\chi}(u) = \det\nolimits_{\mathbf{C}_\ell} \left(1 - F_{*,i} \cdot u |\, \mathrm{H}_i(X, \mathbf{C}_\ell)^{\chi^{-1}}\right).$$

Thus, for all $i$ and $\chi$ as above, we have

$$\begin{aligned} P_i(u) &= \sum_\chi P^{i,\chi}(u) \cdot e_{\chi^{-1}} = \\ &= \sum_\chi \det\nolimits_{\mathbf{C}_p} \left(1 - F^* \cdot u |\, \mathrm{H}_{\mathrm{cris}}^i(X, \mathbf{C}_p)^\chi\right) \cdot e_{\chi^{-1}}. \end{aligned} \quad (\mathrm{a}.12)$$

Let $L$ be a local field of characteristic 0, let $V$ be a finitely generated $L[G]$–module and let $f \in \mathrm{End}_{L[G]}(V)$. We have a (unique up to an $L[f]$–isomorphism) splitting of $V$ as a finite direct sum of $L[G]$–submodules

$$V = \bigoplus_{j \in J} L[f] \cdot v_j, \quad (\mathrm{a}.13)$$

where $L[f] \cdot v_j$ are cyclic $L[f]$–modules, isomorphic to $L[u]/(q_j(u))$, with $q_j(u)$ powers of irreducible polynomials, satisfying $\det(u - f|\, V) = \prod_{j \in J} q_j(u)$. We will denote by $J_0$ the subset of $J$ consisting of those indices $j$ such that the roots of $q_j(u)$ in an algebraic closure $\bar{L}$ of $L$ have valuation 0. Let

$$V_0 = \bigoplus_{j \in J_0} L[f] \cdot v_j, \quad V_1 = \bigoplus_{j \notin J_0} L[f] \cdot v_j.$$

Since $f$ is $G$–invariant, we have a direct sum decomposition of $L[G]$–modules $V = V_0 \oplus V_1$, thus the following divisibility holds in $L[G][u]$:

$$\det\nolimits_{L[G]}(1 - f \cdot u |\, V_0) \mid \det\nolimits_{L[G]}(1 - f \cdot u |\, V). \quad (\mathrm{a}.14)$$

We will apply the general remarks above to a situation involving the étale and crystalline cohomology groups of $X$, described in [3], §3, and which we briefly now summarize. Let $L$ be the quotient field of the Witt vector ring $W(\mathbf{F})$ and let

$$\mathrm{H}_{\mathrm{cris}}^i(X, L) \stackrel{\mathrm{def}}{=} \mathrm{H}_{\mathrm{cris}}^i(X/W(\mathbf{F})) \bigotimes_{W(\mathbf{F})} L,$$

be the crystalline cohomology group of $X$ with coefficients in $L$.

Let $\mathrm{H}^i(X, \mathbf{Q}_p)$ be the étale cohomology groups of $X$ with coefficients in $\mathbf{Q}_p$, defined by $\mathrm{H}^i(X, \mathbf{Q}_p) = \mathrm{H}_i(X, \mathbf{Q}_p)^*$ (functorial dual). Let

$$\mathrm{H}^i(X, \mathbf{C}_p) \stackrel{\mathrm{def}}{=} \mathrm{H}^i(X, \mathbf{Q}_p) \bigotimes_{\mathbf{Q}_p} \mathbf{C}_p.$$

As in the case $\ell \neq p$ explained in detail above, we have $\mathbf{C}_p$–linear isomorphisms

$$\mathrm{H}^i(X, \mathbf{C}_p)^\chi \xrightarrow{\sim} \mathrm{H}_i(X, \mathbf{C}_p)^{\chi^{-1}}, \quad \forall \chi \in \widehat{G},$$

preserving the action of the $q$–power geometric Frobenius morphism, thus

$$\begin{aligned}
P_{i,p}(u) &= \sum_\chi \det\nolimits_{\mathbf{C}_p} \left(1 - F_{*,i} \cdot u \mid \mathrm{H}_i(X, \mathbf{C}_p)^\chi\right) \cdot e_\chi \\
&= \sum_\chi \det\nolimits_{\mathbf{C}_p} \left(1 - F^* \cdot u \mid \mathrm{H}^i(X, \mathbf{C}_p)^\chi\right) \cdot e_{\chi^{-1}}, \forall i = 0, 1, 2.
\end{aligned} \quad (\mathrm{a}.15)$$

There is an exact sequence of $\mathbf{Q}_p[G]$–modules (see [3], Th. 3.2)

$$0 \longrightarrow \mathrm{H}^i(X, \mathbf{Q}_p) \xrightarrow{j} \mathrm{H}^i_{\mathrm{cris}}(X, L) \xrightarrow{1-\phi} \mathrm{H}^i_{\mathrm{cris}}(X, L) \longrightarrow 0,$$

where $\phi$ is the map induced at the cohomology level by the $p$–power geometric Frobenius morphism of $X$ (viewed as an $\mathbf{F}_p$–scheme), and the injective morphism $j : \mathrm{H}^i(X, \mathbf{Q}_p) \longrightarrow \mathrm{H}^i_{\mathrm{cris}}(X, L)$ preserves the action of $F^*$ on both sides (see [3], Lemma 3.3).

Let $\mathrm{H}^i_{\mathrm{cris}}(X, L) = \mathrm{H}^i_{\mathrm{cris}}(X, L)_0 \bigoplus \mathrm{H}^i_{\mathrm{cris}}(X, L)_1$ be the decomposition described above for $V = \mathrm{H}^i_{\mathrm{cris}}(X, L)$ and $f = F^*$. Then the proof of Lemma 3.3 in [3] shows that $\mathrm{H}^i_{\mathrm{cris}}(X, L)_0$ and $\mathrm{H}^i_{\mathrm{cris}}(X, L)_1$ are $\phi$–stable, $j\left(\mathrm{H}^i(X, \mathbf{Q}_p)\right) \subseteq \mathrm{H}^i_{\mathrm{cris}}(X, L)_0$, and via $j$ we have an equality of $L[G]$–modules

$$\mathrm{H}^i_{\mathrm{cris}}(X, L)_0 = \mathrm{H}^i(X, \mathbf{Q}_p) \bigotimes_{\mathbf{Q}_p} L.$$

Obviously this last equality can be written as

$$\mathrm{H}^i_{\mathrm{cris}}(X, L)_0 = \mathrm{H}^i(X, \mathbf{Q}_p) \bigotimes_{\mathbf{Q}_p[G]} L[G].$$

This shows that

$$\det\nolimits_{\mathbf{Q}_p[G]} \left(1 - F^* \cdot u \mid \mathrm{H}^i(X, \mathbf{Q}_p)\right) = \det\nolimits_{L[G]} \left(1 - F^* \cdot u \mid \mathrm{H}^i_{\mathrm{cris}}(X, L)_0\right),$$

which, according to (a.14), implies that, inside $\mathbf{C}_p[G][u]$, we have

$$\det\nolimits_{\mathbf{C}_p[G]} \left(1 - F^* \cdot u \mid \mathrm{H}^i(X, \mathbf{C}_p)\right) \mid \det\nolimits_{\mathbf{C}_p[G]} \left(1 - F^* \cdot u \mid \mathrm{H}^i_{\mathrm{cris}}(X, \mathbf{C}_p)\right).$$

If we break the last divisibility into $\chi$–components, we obtain the following divisibilities in $\mathbf{C}_p[u]$:

$$\det_{\mathbf{C}_p}\left(1 - F^* \cdot u|\, \mathrm{H}^i\left(X, \mathbf{C}_p\right)^\chi\right) \mid \det_{\mathbf{C}_p}\left(1 - F^* \cdot u|\, \mathrm{H}^i_{\mathrm{cris}}\left(X, \mathbf{C}_p\right)^\chi\right), \quad \forall \chi \in \widehat{G}.$$

According to (a.12) and (a.15), these relations imply that there exist polynomials $Q_i(u) \in \mathbf{C}_p[G][u]$ such that

$$P_i(u) = Q_i(u) \cdot P_{i,p}(u), \quad \forall i = 0, 1, 2.$$

In order to prove statements (3) and (4) in our theorem, we will just need to show that $Q_i(u) \in \mathbf{Z}_p[1/g][G][u]$ and that, in the case when $K/k$ is a constant field extension, $Q_i(u) \in \mathbf{Z}_p[G][u]$. These statements will easily follow from the next elementary lemma.

**Lemma A.5.** *Let $S \subseteq R$ be two commutative rings with $1$. Let $u$ be a variable and let $f(u) \in S[u]$, $g(u) \in S[u]$ and $h(u) \in R[u]$, such that $g(0) = 1$, and $f(u) = g(u) \cdot h(u)$. Then $h(u) \in S[u]$ as well.*

**Proof.** Let $f(u) = \sum_{i \geq 0} f_i u^i$, $g(u) = \sum_{i \geq 0} g_i u^i$, and $h(u) = \sum_{i \geq 0} h_i u^i$, with $f_i, g_i \in S$, $h_i \in R$ and $g_0 = 1$. The hypotheses of our lemma imply that $h_0 = f_0$ and $f_1 = g_1 h_0 + h_1$. Since $f_0, f_1, g_1 \in S$, these equalities imply that $h_0, h_1 \in S$. Let us supose that $h_0, h_1, \ldots, h_s \in S$, for $0 \leq s \leq \deg(h) - 1$. We have a relation

$$f_{s+1} = h_{s+1} + (h_s g_1 + h_{s-1} g_2 + \ldots),$$

which shows that

$$h_{s+1} = f_{s+1} - (h_s g_1 + h_{s-1} g_2 + \ldots).$$

The induction hypothesis now shows that $h_{s+1} \in S$.  □

We are now ready to conclude the proof of Theorem 1.7.4.1 (3) and (4).

**Proof of (3).** According to (1) and (2) in our Theorem, the desired statement follows from Lemma A.5 applied to $f(u) = P_i(u)$, $g(u) = P_{i,p}(u)$, $h(u) = Q_i(u)$, $S = \mathbf{Z}_p[1/g][G][u]$ and $R = \mathbf{C}_p[G][u]$.  □

**Proof of (4).** Let $K/k$ be a finite, constant field extension of function fields. Remark 1 above shows that

$$\Theta(u) \in \mathbf{Z}[G][[u]]. \tag{a.16}$$

The proof of Lemma 4.2.1 shows that

$$P_0(u) = 1 - \sigma^{-1} u \in \mathbf{Z}[G][u], \quad P_2(u) = 1 - q\sigma^{-1} u \in \mathbf{Z}[G][u],$$

where $\sigma$ is the distinguished generator of $G$ described in §4.1. According to (a.7) and (a.16), we therefore have

$$P_i(u) \in \mathbf{Z}[G][u], \quad \forall i = 0, 1, 2. \tag{a.17}$$

Since the $\mathrm{H}_i(X, \mathbf{Z}_p)$'s are free $\mathbf{Z}_p[G]$–modules in the constant field extension case (see Lemma 4.1.1), we also have

$$P_{i,p}(u) = \det\nolimits_{\mathbf{Z}_p[G]}(1 - F_* \cdot u | \mathrm{H}_i(X, \mathbf{Z}_p)) \in \mathbf{Z}_p[G][u], \quad \forall i = 0, 1, 2. \tag{a.18}$$

According to (a.17) and (a.18), statement (4) in our Theorem follows from Lemma A.5 applied to $f(u) = P_i(u)$, $g(u) = P_{i,p}(u)$, $h(u) = Q_i(u)$, $S = \mathbf{Z}_p[G][u]$ and $R = \mathbf{C}_p[G][u]$.  $\square$

MATHEMATICAL SCIENCES RESEARCH INSTITUTE, 1000 CENTENNIAL DRIVE, BERKELEY, CA, 94720–5070, USA
   *E-mail address*: : popescu@msri.org